\documentclass[a4paper,11pt]{article}

\usepackage{geometry,enumerate,amsmath,amssymb,latexsym,citesort}
\usepackage[all]{xy}
\newtheorem{example}{Example}[section]
\newtheorem{definition}[example]{Definition}
\newtheorem{proposition}[example]{Proposition}
\newtheorem{theorem}[example]{Theorem}
\newtheorem{remark}[example]{Remark}
\newtheorem{corollary}[example]{Corollary}

\newenvironment{proof}{\noindent {\bf Proof }}{\hfill $\Box$

\mbox{}}

\def\rho{\varrho}

\def\del{\partial}
\def\A{\alpha}
\def\B{\beta}
\def\eps{\varepsilon}
\def\epsilon{\varepsilon}

\def\half{\frac{1}{2}}

\def\le{\leqslant}
\def\ge{\geqslant}

\def\d{\mathbf{d}}

\def\del{\partial}
\def\A{\alpha}
\def\B{\beta}
\def\eps{\varepsilon}
\def\epsilon{\varepsilon}

\def\half{\frac{1}{2}}

\def\geq{\geqslant}
\def\leq{\leqslant}
\def\ge{\geqslant}
\def\le{\leqslant}

\def\dip{_{i+1}}
\newcommand{\directs}[2]{\def\objectstyle{\scriptstyle}  \objectmargin={0pt}
\xy
(0,4)*+{}="a",(0,-2)*+{\rule{0em}{1.5ex}#2}="b",(7,4)*+{\;#1}="c"
\ar@{->} "a";"b" \ar @{->}"a";"c" \endxy }

\newcommand{\xdirects}[2]{\def\objectstyle{\scriptstyle}  \objectmargin={0pt}
\xy
(0,0)*+{}="a",(0,-6)*+{\rule{0em}{1.5ex}#2}="b",(7,0)*+{\;#1}="c"
\ar@{->} "a";"b" \ar @{->}"a";"c" \endxy }

 \def\c{\mathbin{\#}}
 
 \def\d{\partial}
 \def\Hom{\mathop{\rm Hom}\nolimits}
 \def\id{\mathop{\rm id}\nolimits}
 \def\Im{\mathop{\rm Im}\nolimits}

 \def\epsilon{\varepsilon}
 \def\cat{\mathsf{Cat}}
  \def\Cat{\mathsf{CAT}}
  \def\glob{\bigcirc}
  \def\dash{\text{-}}
\def\omcat{cubical $\omega$-categor}
\def\ocat{$\omega$-$\mathsf{Cat}^\Box$}
\def\oCat{$\omega$-$\mathsf{CAT}^\Box$}
\def\ocatm{\omega\mbox{-}\mathsf{Cat}^\Box}
\def\oCatm{\omega\mbox{-}\mathsf{CAT}^\Box}

\def\eps{\varepsilon}
\def\cub{\mathsf{Cub}}
\def\Cub{\mathsf{CUB}}
\def\bI{\mathbb{I}}

\def\P{\mathsf{P}}
\begin{document}

\title{ Multiple categories: \\ the equivalence of a  globular
and a cubical approach}
\author{Fahd Ali Al-Agl\footnote{We have been unable to contact Dr.
Al-Agl and so the responsibility for the paper rests with the last
two authors.}, \\Um-Alqura University,\\ Makkah\\Saudi Arabia \and Ronald Brown, \\ School of Informatics, \\ Mathematics Division, \\ University of
Wales,\\ Bangor, Gwynedd LL57 1UT, \\    United Kingdom.\\ email:
r.brown@bangor.ac.uk \and Richard Steiner, \\ Department of Mathematics, \\
University of Glasgow, \\University Gardens, \\ Glasgow G12 8QW
\\ United Kingdom\\ email: r.steiner@maths.gla.ac.uk
 }
\date{17 December, 2001}
  \maketitle
\begin{center}
Accepted for Advances in Mathematics.

{\em
Dedicated to Professor Philip J. Higgins for his 75th birthday}
\end{center}

\begin{abstract}
We show the equivalence of two kinds of strict multiple category,
namely the well known globular $\omega$-categories, and the
cubical $\omega$-categories with connections.
\end{abstract}


\tableofcontents

\section*{Introduction}

An essential feature for the possibility of `higher dimensional
group theory' (see the expository article \cite{B2}) is the
extension of the domain of discourse from groups to groupoids,
that is from a set with a binary operation defined on all
elements, to a set with an operation defined only on pairs
satisfying a geometric condition. This fact itself leads to
various equivalent candidates for `higher dimensional groups',
namely those based on different geometric structures, for example
balls, globes, simplices, cubes and even polyhedra.  The proofs of
these equivalences are  non trivial -- the basic intuitions derive
from the foundations of relative homotopy theory. Some of these
equivalences have proved crucial for the applications: theorems
may be easily proved in one context and then transferred into
another, more computational context. Notable examples are the
advantages of cubical methods for providing both a convenient
`algebraic inverse to subdivision', for use in local-to-global
problems \cite{BH2a}, and also a simple monoidal closed structure,
which may then be translated into other situations \cite{BHtens}.

It has proved important to extend these ideas from groupoids to
categories.  The standard notion of (strict) higher dimensional
category is that of globular $\omega$-category. Our main result is
that there is an adjoint equivalence of categories $$
\lambda:\text{globular } \omega\text{-categories}
 \xymatrix{\ar@<0.5ex>[r]&\ar@<0.5ex>[l]}
  \text{cubical }\omega\text{-categories with connections }
: \gamma .
$$ Precise definitions are given below. The proof has interest because
it is certainly much harder than the groupoid case, and because at
one stage it uses braid relations among some key basic folding
operations (Proposition \ref{5.1}, Theorem \ref{5.2}). The
equivalence between the two forms  should prove useful. In
section 9 we use this equivalence to define the notion of
`commutative $n$-cube'.  In section 10 we follow methods of Brown
and Higgins in \cite{BHtens} to show that cubical
$\omega$-categories with connection form a monoidal closed
category. The equivalence of categories transfers this structure
to the globular case -- the resulting internal hom in the
globular case gives various higher dimensional forms of `lax
natural transformation'.  Cubical $\omega$-categories with
connection have been applied to concurrency theory by E. Goubault
\cite{Goub} and by P. Gaucher \cite{Gauch}, and again relations
with the globular case are important for these studies.

The origin of this equivalence is as follows.

In developing the algebra of double groupoids as a framework for
potential 2-dimensional Van Kampen Theorems, Brown and Spencer in
\cite{Br-Spe} were led to the notion of double groupoid with an
extra structure of `connection' -- this was essential to obtain an
equivalence of such a double groupoid with the classical notion of
crossed module. This structure was also essential for the proof of
the 2-dimensional Van Kampen Theorem given by Brown and Higgins in
\cite{BH2d}.

The double groupoid case was generalised by Brown and Higgins
\cite{BH1,BH2}  to give an equivalence between crossed complexes
and what were called there `$\omega$-groupoids', and which we here
call `cubical $\omega$-groupoids with connections'. It was also
proved in \cite{BH3} that crossed complexes are equivalent to what
were there called `$\infty$-groupoids', and which we here call
`globular  $\omega$-groupoids', following current fashions. Thus
the globular and cubical cases of $\omega$-groupoids were known in
1981 to be equivalent, but the proof was via the category of
crossed complexes.

Other equivalences with crossed complexes were established, for
example with: cubical $T$-complexes \cite{BH1,BH4};  simplicial
$T$-complexes by Nick Ashley \cite{Ash}; and polyhedral
$T$-complexes by David Jones \cite{J}. In $T$-complexes the basic
concept is taken to be that of {\em thin elements} which determine
a strengthening of the Kan extension condition. The notion of
simplicial $T$-complex is due to Keith Dakin \cite{Da}.

Spencer observed in \cite{Spe} that the methods of \cite{Br-Spe}
allowed an equivalence between 2-categories and double categories
with connections, using an `up-square' construction of Bastiani
and Ehresmann \cite{BE,E}, but he gave no details. The full
details of this have been recently given by Brown and Mosa in
\cite{Br-Mo}.

The thesis of Mosa in 1987 \cite{Mo} attempted to give an
equivalence between crossed complexes of algebroids and cubical
$\omega$-algebroids, and while this was completed in dimension 2
even the case of dimension 3 proved hard, though some basic
methods were established.

This result raised the question of an equivalence between the
globular  $\omega$-categories defined in 1981 in \cite{BH3} and an
appropriate form of cubical $\omega$-categories with connections,
of which a definition was fairly easy to formulate as an extension
of the previous definition of cubical $\omega$-groupoid. This
problem was taken up in Al-Agl's thesis of 1989 \cite{AA}. The
central idea, based on the groupoid methods of \cite{BH2}, was to
define a `folding operation' $\Phi$ from a cubical
$\omega$-category $G$ to the globular $\omega$-category $\gamma G$
it contained. This definition was successfully accomplished, but
the problem of establishing some major properties of $\Phi$, in
particular the relation with the category structures, was solved
only up to dimension 3. That is, the conjectured equivalence was
proved in dimension 3.

Steiner pursued the work of Al-Agl, and their joint paper
\cite{AAS} does prove that globular $\omega$-categories are
equivalent to cubical sets with extra structure, but, as stated in
that paper, this extra structure is not described in finitary
terms. Later, Steiner was stimulated by renewed interest in the
cubical case coming from concurrency theory, in work of E.
Goubault \cite{Goub} and P. Gaucher \cite{Gauch}, and by the
publication of the 2-dimensional case by Brown and Mosa in
\cite{Br-Mo}.  He completed the programme given in \cite{AA} and
informed Brown, who announced the result at the Aalborg `Workshop
on Geometric and Topological Methods in Concurrency' in June,
1999. This paper is the result. It proves the conjecture implicit
in \cite{AA}, that a globular $\omega$-category is equivalent to a
cubical set with extra structure directly analogous to the
structure for cubical $\omega$-groupoids given in \cite{BH1,BH2}.

There is considerable independent work on globular
$\omega$-categories. The thesis of Sjoerd Crans \cite{Cr} already
contains the adjoint pair $(\lambda,\gamma)$ and also the closed
monoidal structure on  the category of globular
$\omega$-categories. It also seems to be the first time that that
the cube category (without connections) together with its
$\omega$-category realisation is explicitly defined by generators
and relations.

The work  in Australia by Ross Street
\cite{Street1,Street2,Street2a,Street3} has an initial aim to
determine a simplicial nerve $NX$ of a globular $\omega$-category
$X$. This developed into finding extra structure on $NX$ so that
$N$ gave an equivalence  between $\omega$-categories and certain
structured simplicial sets, analogous to Ashley's equivalence
\cite{Ash} between $\omega$-groupoids and simplicial
$T$-complexes. It is stated in \cite{Street3} that this programme
has been completed by Dominic Verity, to verify the conjecture
stated in \cite{Street2a}. Street tells us that Verity also knew
the equivalence proved in the present paper, but we have no
further information. We also mention that Street's paper
\cite{Street2a} {\em implicitly} contains our basic  proposition
(3.2), namely that the cells of the $n$-categorical $n$-cube
compose in such a way that they give rise to the hemispherical
(i.e. globular) decomposition $\partial^\pm_1 \Phi_n$ of the
$n$-cube.

\section{$\omega$-categories}

An $\omega$-category \cite{BH3,Street1,Ste1} arises when a
sequence of categories $C_0,C_1,\ldots\,$ all have the same set of
morphisms~$X$, the various category structures commute with one
another, the identities for~$C_p$ are also identities for~$C_q$
when $q>p$, and every member of~$X$ is an identity for some~$C_p$.
We write~$\c_p$ for the composition in~$C_p$. Given $x\in X$, we
write $d^-_p x$ and $d^+_p x$ for the identities of the source and
target of~$x$ in~$C_p$, so that $d^-_p x\c_p x=x\c_p d^+_p x=x$.
The structure can be expressed in terms of $X$, $\c_p$ and
the~$d^\alpha_p$ as follows.

\begin{definition}\label{1.1} {\em An {$\omega$-category\/} is a set~$X$
together with unary operations $d^-_p$,~$d^+_p$ and partially
defined binary operations~$\c_p$ for $p=0$,~$1$,~\dots\ such that
the following conditions hold:
\begin{enumerate}[(i)]
  \item  $x\c_p y$ is defined if and only if $d^+_p x=d^-_p y$;
  \item  $d^\beta_q d^\alpha_p x=\begin{cases}d^\beta_q x&
  \text{for } q<p,\\
    d^\alpha_p x& \text{for }q\geq   p; \end{cases} $
    \item   if $x\c_p y$ is defined then
 \begin{align*}
 d^-_p(x\c_p y)& =d^-_p x,\\ d^+_p(x\c_p y)& =d^+_p y,\\
 d^\beta_q(x\c_p y)& =d^\beta_q x\c_p d^\beta_q y\ \text{for }
 q\neq p;
 \end{align*}
 \item   $d^-_p x\c_p x=x\c_p d^+_p x=x$;
 \item   $(x\c_p y)\c_p z=x\c_p(y\c_p z)$ if either side is defined;
 \item if $p\neq q$ then
$$(x\c_p y)\c_q(x'\c_p y')=(x\c_q x')\c_p(y\c_q y')$$
whenever both sides are defined;
\item  for each $x\in X$ there is a {\it dimension\/} $\dim x$ such
that $d^\alpha_p x=x$ if and only if $p\geq\dim x$.
\end{enumerate}
               }\end{definition}

\begin{definition} \label{1.2} {\em An {\it $\omega$-category of sets\/}
is an $\omega$-category~$X$ whose members are sets such that
$x\c_p y=x\cup y$ whenever $x\c_p y$ is defined in~$X$.}
\end{definition}

The theory of pasting in $\omega$-categories \cite{Ste1,Street2a}
associates $\omega$-categories of sets $M(K)$ with simple
presentations to certain complexes~$K$; the members of $M(K)$ are
subcomplexes of~$K$. Various types of complexes have been
considered, but they certainly include the cartesian products of
directed paths, and we will now describe the theory in that case.

Let $n$ be a non-negative integer. We represent a directed path of
length~$n$ by the closed interval $[0,n]$; the vertices are the
singleton subsets $\{0\}$, $\{1\}$, \dots, $\{n\}$ and the edges
are the intervals $[0,1]$, $[1,2]$, \dots, $[n-1,n]$, where
$[m-1,m]$ is directed from $m-1$ to~$m$. We write
$$d^-[m-1,m]=\{m-1\},\ d^+[m-1,m]=\{m\}.$$

Now let $K=K_1\times\cdots\times K_p$ be a cartesian product of
directed paths. A product $\sigma=\sigma_1\times\cdots\times
\sigma_p$, where $\sigma_i$~is a vertex or edge in~$K_i$, is
called a {\it cell\/} in~$K$. We can write a cell~$\sigma$ in the
form
 $$\sigma=P_0\times e_1\times P_1\times e_2\times
 P_2\times\cdots\times  P_{q-1}\times e_q\times P_q, $$
where the~$P_j$ are products of vertices and the~$e_j$ are edges;
the {\it dimension\/} of~$\sigma$ is then~$q$. The codimension~$1$
faces of~$\sigma$ are the subsets got by replacing one edge
factor~$e_j$ with $d^-e_j$ or $d^+e_j$. The faces with $d^-e_1$ or
$d^+e_2$ or $d^-e_3$ or \dots\ are called {\it negative}, and the
faces with $d^+e_1$ or $d^-e_2$ or $d^+e_3$ or \dots\ are called {\it
positive}. The theory of pasting gives us the following result.

\begin{theorem} \label{1.3} Let $K$ be a cartesian product of directed
paths. Then there is an $\omega$-category $M(K)$ of subsets of~$K$
with the following presentation\/{\rm:} the generators are the
cells of~$K${\rm;} if $\sigma$ is a cell of dimension~$q$ then
there are relations $d^-_q\sigma=d^+_q\sigma=\sigma${\rm;} if
$\sigma$ is a cell of dimension~$q$ with $q>0$ then there are
relations saying that $d^-_{q-1}\sigma$ and $d^+_{q-1}\sigma$ are
the unions of the negative and positive faces of~$\sigma$
respectively. Every member of $M(K)$ is an iterated composite of
cells.
\end{theorem}

We will now describe the main examples.

\begin{example} \label{1.4}{\em  We write $I=[0,1]$ and
$I^n=[0,1]^n$ for $n\geq 1$; for completeness we also write
$I^0=[0,0]$. In this notation $M(I^0)=\{I^0\}$ and $M(I)=\{I,d^-_0
I,d^+_0 I\}$; there are no members other than the generating
cells. There are morphisms
 $$\check\d^-,\check\d^+\colon M(I^0)\to
 M(I),\quad \check\epsilon\colon M(I)\to M(I^0)$$
given by
\begin{gather*}
 \check\d^\alpha(I^0)=d^\alpha_0 I,
 \\  \check\epsilon(I)=\check\epsilon(d^\alpha_0 I)=I^0.
 \end{gather*}
}\end{example}

\begin{example} \label{1.5} {\em The members of $M([0,2])$
are the cells and the composite $$[0,2]=[0,1]\c_0[1,2].$$ There
are morphisms
 $$\check\iota^-,\check\iota^+,\check\mu\colon
 M(I)\to M([0,2])$$
given by
\begin{gather*}
 \check\iota^-(d^-_0 I)=\{0\},\
 \check\iota^-(I)=[0,1],\
 \check\iota^-(d^+_0 I)=\{1\},\\  \check\iota^+(d^-_0 I)=\{1\},\
 \check\iota^+(I)=[1,2],\
 \check\iota^+(d^+_0 I)=\{2\},\end{gather*}
and
 $$\check\mu(d^-_0 I)=\{0\},\
 \check\mu(I)=[0,2],\
 \check\mu(d^+_0 I)=\{2\}.$$     }
\end{example}

\begin{example} \label{1.6} {\em The members of $M(I^2)$ are
the cells and the composites
 $$d^-_1 I^2=(d^-_0 I\times I)\c_0(I\times d^+_0 I),\quad
 d^+_1 I^2=(I\times d^-_0 I)\c_0(d^+_0 I\times I).$$
There are morphisms $\check\Gamma^+,\check\Gamma^-\colon M(I^2)\to
M(I)$ given by
\begin{gather*}
 \check\Gamma^{\alpha}(d^{-\alpha}_0 I\times d^{-\alpha}_0 I)
 =\check\Gamma^{\alpha}(d^{-\alpha}_0 I\times I)
 =\check\Gamma^{\alpha}(d^{-\alpha}_0 I\times d^{\alpha}_0 I)
 =\check\Gamma^{\alpha}(I\times d^{-\alpha}_0 I) \\
 =\check\Gamma^{\alpha}(d^{\alpha}_0 I\times d^{-\alpha}_0 I)
 =d^{-\alpha}_0 I,\\  \check\Gamma^{\alpha}(I^2)
 =\check\Gamma^{\alpha}(I\times d^{\alpha}_0 I)
 =\check\Gamma^{\alpha}(d^{\alpha}_0 I\times I)
 =\check\Gamma^{\alpha}(d^-_1 I^2)
 =\check\Gamma^{\alpha}(d^+_1 I^2)
 =I,\\  \check\Gamma^{\alpha}(d^{\alpha}_0 I\times d^{\alpha}_0 I)
 =d^{\alpha}_0 I.\end{gather*}  }
\end{example}

For cartesian products of members of the $\omega$-categories that
we are considering, we have the following result.

\begin{theorem} \label{1.7} Let $K$~and~$L$ be cartesian products of
directed paths, let $x$ be a member of $M(K)$, and let $y$ be a
member of $M(L)$. Then $x\times y$ is a member of $M(K\times L)$
and
$$d^\alpha_p(x\times y)
 =\bigcup_{i=0}^p
 \left(d^\alpha_i x\times d^{(-)^i\alpha}_{p-i}y\right).$$
\end{theorem}

This has the following consequence.

\begin{theorem} \label{1.8} {\rm (i)} Let $K$, $K'$, $L$, $L'$ be
cartesian products of directed graphs, and let $f\colon M(K)\to
M(K')$ and $g\colon M(L)\to M(L')$ be morphisms of
$\omega$-categories. Then there is a unique morphism $$f\otimes
g\colon M(K\times L)\to M(K'\times L')$$ such that $$(f\otimes
g)(x\times y)=f(x)\times g(y)$$ for $x\in M(K)$ and $y\in M(L)$.

{\rm (ii)} The assignments
 $$\bigl(M(K),M(L)\bigr)\mapsto M(K\times L),\quad
 (f,g)\mapsto f\otimes g$$
form a bifunctor.
\end{theorem}

\begin{proof}(i) From the presentation of $M(K\times L)$ and
Theorem \ref{1.7}, there is a unique morphism $f\otimes g$ such that
$(f\otimes g)(x\times y)=f(x)\times g(y)$ when $x$~and~$y$ are
cells. The formula then holds for a general product $x\times y$
because it is a composite of cells.

(ii) One can check bifunctoriality by considering the values of
the appropriate morphisms on generators.
\end{proof}

By applying the tensor product construction, we obtain further
morphisms.

\begin{example} \label{1.9} {\em Let $\id^r$ denote the
identity morphism from $M(I^r)$ to itself. There are morphisms
$$\check\d^-_i,\check\d^+_i\colon M(I^{n-1})\to M(I^n)\quad
 (1\leq i\leq n)$$
given by
$$\check\d^\alpha_i=\id^{i-1}\otimes\check\d^\alpha\otimes\id^{n-i};$$
there are morphisms
$$\check\epsilon_i\colon M(I^n)\to M(I^{n-1})\quad
 (1\leq i\leq n)$$
given by
$$\check\epsilon_i=\id^{i-1}\otimes\check\epsilon\otimes\id^{n-i};$$
there are morphisms
$$\check\iota^-_i,\check\iota^+_i,\check\mu_i\colon
 M(I^n)\to M(I^{i-1}\times[0,2]\times I^{n-i})\quad
 (1\leq i\leq n)$$
given by
$$
 \check\iota^\alpha_i=\id^{i-1}\otimes\check\iota^\alpha\otimes\id^{n-i},
 \quad
 \check\mu_i=\id^{i-1}\otimes\check\mu\otimes\id^{n-i};$$
there are morphisms
$$\check\Gamma^+_i,\check\Gamma^-_i\colon M(I^n)\to M(I^{n-1})\quad
 (1\leq i\leq n-1)$$
given by
$$\check\Gamma^{\alpha}_i
 =\id^{i-1}\otimes\check\Gamma^{\alpha}\otimes\id^{n-i-1}.$$ }
\end{example}

Most of the morphisms in Example \ref{1.9} map generators to
generators, and one can verify their existence directly from
Theorem \ref{1.3}. The exceptions are the~$\check\mu_i$, for which
Theorem \ref{1.8} is really necessary.

\begin{remark} \label{1.10}{\em  Suppose that $K$ is an $n$-dimensional
product of directed paths. Then $K$ can be got from a family of
$n$-cubes by gluing along $(n-1)$-dimensional faces. From the
presentation of $M(K)$, one sees that it is the colimit of a
corresponding diagram in which the morphisms have the form
$\check\d^\alpha_i\colon M(I^{n-1})\to M(I^n)$. In particular,
$\check\iota^-$~and~$\check\iota^+$ exhibit $M([0,2])$ as the
push-out of $$\xymatrix@1{
 M(I) & M(I^0)\ar[l]_-{\check\d^+}\ar[r]^-{\check\d^-} &
 M(I),}$$
and $\check\iota^-_i$~and~$\check\iota^+_i$ exhibit
$M(I^{i-1}\times[0,2]\times I^{n-i})$ as the push-out of
$$\xymatrix@1{
 M(I^n) & M(I^{n-1})\ar[l]_-{\check\d^+_i}\ar[r]^-{\check\d^-_i} &
 M(I^n).}$$                  }
\end{remark}

\section{Cubical $\omega$-categories with connections}

Suppose that $X$~is an $\omega$-category. There is then a sequence
of sets $$(\lambda X)_n=\Hom[M(I^n),X]\quad (n=0,1,\ldots\,),$$
and the morphisms of Example \ref{1.9} induce functions between
the $(\lambda X)_n$. It turns out that the $(\lambda X)_n$ form a
cubical $\omega$-category with connections in the sense of the following
definition. This definition is found in \cite{AA}. The origin is
in the definition of what was called `$\omega$-groupoid' in
\cite{BH1,BH2}, where the justification was the equivalence with
crossed complexes ({\em loc. cit.}) and the use in the formulation
and proof of a generalised Van Kampen Theorem \cite{BH1a,BH2a}.
The corresponding definition for categories arose out of the work
of Spencer \cite{Spe} and of Mosa \cite{Mo}.

Let  $K$  be a cubical set, that is, a family of  sets  $\{ K_n;n
\ge 0\} $ with for $n \ge 1 $ face maps  $ \del _i^\A:K_n  \to
K_{n-1} \; (i = 1,2,\ldots,n;\, \alpha = +,-)$  and degeneracy
maps $\epsilon_i:K_{n-1}     \to   K_n \; (i = 1,2,\ldots,n)$
satisfying the usual cubical relations:

\begin{alignat*}{2}
         \del_i^\A \del_j^\B &=  \del_{j-1}^\B \del_i^{\alpha}
        &&\hspace{-5em}(i<j),
        \tag*{(2.1)(i)} \\
        \epsilon_i \epsilon _j &=  \epsilon _{j+1} \epsilon _i && \hspace{-5em}(i \le j),
          \tag*{(2.1)(ii)} \\
\del ^{\alpha}_i \epsilon _j &=
                  {\begin{cases} \eps_{j-1} \del _i^\A &\hspace{8em} (i<j)  \\
                          \eps_{j} \del _{i-1}^\A & \hspace{8em}(i>j)  \\
                             \mathrm{id} &\hspace{8em} (i=j)
                  \end{cases}} && \tag*{(2.1)(iii)}   \\
\intertext{ We say that  $K$  is a  {\em cubical  set with
connections}  if  for $n \ge 0$ it  has additional structure maps
(called {\em connections}) $\Gamma_i^+,\Gamma_i^- :K_{n} \to
K_{n+1} \; (i = 1,2,\ldots,n)$ satisfying the relations:}
       \Gamma_i^\A   \Gamma_j^\B  & =  \Gamma_{j+1}^\B\Gamma_i^\A &&
        \hspace{-5em}  (i  <  j)
       \tag*{(2.2)(i)} \\
         \Gamma_i^\A   \Gamma_i^\A  & =  \Gamma_{i+1}^\A\Gamma_i^\A &&
       \tag*{(2.2)(ii)} \\
       \Gamma_i^\A   \epsilon_j & = {\begin{cases}  \epsilon_{j+1}\Gamma_i^\A &\hspace{8em} (i < j)\\
                                     \epsilon_{j}\Gamma_{i-1}^\A &\hspace{8em}(i > j)
                                  \end{cases}}&& \tag*{(2.2)(iii)} \\
        \Gamma_j^\A \eps_j &= \eps^2_j=\eps_{j+1}\eps_j, &&  \tag*{(2.2)(iv)} \\
        \del^\A _i \Gamma_j^\B &= {\begin{cases} \Gamma_{j-1}^\B\del^\A_i
        & \hspace{8em}(i<j) \\
                              \Gamma_{j}^\B\del^\A_{i-1}  &\hspace{8em} (i> j+1),
                              \end{cases}}&& \tag*{(2.2)(v) }\\
        \del^\A_j\Gamma_j^\A&= \del_{j+1}^\A   \Gamma _j^\A  =  id, && \tag*{(2.2)(vi)} \\
        \del^\A_j \Gamma^{-\A} _j&=
        \del^\A_{j+1}  \Gamma^{-\A} _j  = \epsilon_j  \del^\A_j.
        &&        \tag*{(2.2)(vii)}
\end{alignat*}
The connections are to be thought  of  as  extra  `degeneracies'.
(A degenerate cube of type  $ \epsilon_j  x$  has a pair of
opposite  faces  equal and all other faces degenerate.  A cube of
type  $ \Gamma_i^\A  x$  has a pair  of adjacent faces equal and
all other faces of type $\Gamma_j^\A  y$  or $\epsilon_j y$ .)
Cubical complexes with this, and other, structures  have  also
been considered by Evrard \cite{Ev}.

The prime example of  a  cubical   set  with  connections  is the
singular cubical complex  $KX$  of a space   $X$.   Here for $n
\ge 0$ $K_n$ is  the  set  of singular $n$-cubes in $X$ (i.e.
continuous maps $I^n   \to   X$)  and the connection $
\Gamma_i^\A :K_{n } \to   K_{n+1}$  is induced by the map
$\gamma_i^\A   : I^{n+1} \to
 I^{n}$    defined by
   $$   \gamma _i^\A (t_1 ,t_2 ,\ldots,t_{n+1} ) =
             (t_1 ,t_2 ,\ldots,t_{i-1},A(t_i ,t_{i+1}),t_{i+2},\ldots,t_{n+1} )
             $$
 where $A(s,t)=\max(s,t), \min(s,t)$ as $\A=-,+$ respectively.
 Here are pictures of $\gamma^\alpha_1 : I^2 \to I^1$ where the
 internal lines show lines of constancy of the map on $I^2$.
\vspace{1.5in}
 \begin{center}
\setlength{\unitlength}{0.1in}
\begin{picture}(0,-20)(0,-10)
\put(-5,2){\makebox{$\gamma^{-}_1 = $}}
\put(-17,2){\makebox{$\gamma^{+}_1 = $}}

\put(0,0){\framebox(4,4){}} \put(0,1){\line(2,0){3}}
\put(0,2){\line(3,0){2}} \put(2,4){\line(0,-2){2}}
\put(0,3){\line(3,0){1}} \put(2,4){\line(0,-2){2}}
\put(1,4){\line(0,-2){1}} \put(3,4){\line(0,-2){3 }}
\put(5,2){\directs{1}{2}}

\put(-12,0){\framebox(4,4){}} \put(-11,3){\line(2,0){3}}
\put(-11,3){\line(0,-1){3}} 
\put(-10,2){\line(3,0){2}} \put(-10,2){\line(0,-2){2}}
\put(-9,1){\line(0,-2){1}} \put(-9,1){\line(1,0){1 }}
\end{picture}
\end{center}

The complex  $KX$  has some further relevant structure, namely the
composition of $n$-cubes in the $n$  different directions.
Accordingly, we define a {\it cubical complex with connections
and compositions} to be a cubical set  $K$  with connections in
which each $K_n$ has  $n$ partial compositions $\circ _j\; (j =
1,2,\ldots,n)$ satisfying the following axioms.

If  $a,b \in K_n$, then  $a\circ _j b$  is defined if and only if    $\del^-_j   b =
\del^+_j  a$  , and then
\begin{equation} \begin{cases} \del^-_j  (a\circ _j b) = \del^-_ja & \\
                 \del^+_j  (a\circ _j b) = \del^+_jb & \end{cases}
                 \qquad
 \del^\A_i  (a\circ _j b) =  \begin{cases} \del^\A_ja\circ _{j-1}\del^\A_i b &(i<j) \\
                 \del^\A_i  a\circ _j \del^\A_i b& (i>j), \end{cases}
                 \tag*{(2.3)}    \end{equation}

 {\em The  interchange laws}.  If  $i \ne j$  then
\begin{equation}
      (a\circ _i b) \circ _j  (c\circ _i d) = (a\circ _j c) \circ _i  (b\circ _j d)
      \tag*{(2.4)}
\end{equation}
whenever both sides are defined. (The diagram
$$
  \begin{bmatrix}
               a &b \\c&d
  \end{bmatrix}     \quad \directs{i}{j}
$$ will be used to indicate that both sides of the above equation
are  defined and also to denote the unique composite of the four
elements.)

     If  $i \ne j$  then
\begin{align*}
           \epsilon_i(a\circ _j b) &= \begin{cases}
           \eps_ia \circ _{j+1} \eps_ib & (i \le j) \\
           \eps_ia \circ _j\eps_ib & (i >j) \end{cases} \tag*{(2.5)} \\
       \Gamma^\A _i (a\circ _j b)& =  \begin{cases}
           \Gamma^\A_ia \circ _{j+1} \Gamma^\A_ib & (i < j) \\
           \Gamma^\A_ia \circ _j\Gamma^\A_ib & (i >j) \end{cases}
           \tag*{(2.6)(i)} \\
       \Gamma^+_j(a\circ _jb)&=  \begin{bmatrix}\Gamma^+_ja & \eps_j a\\
       \eps_{j+1} a & \Gamma^+_j b \end{bmatrix} \quad \directs{j}{j+1}  \tag*{(2.6)(ii)}\\
       \Gamma^-_j(a\circ _jb)&=  \begin{bmatrix}\Gamma^-_ja & \eps_{j+1} b\\
       \eps_{j} b & \Gamma^-_j b \end{bmatrix}\quad \directs{j}{j+1} \tag*{(2.6)(iii)} \\
\end{align*}        These last two equations are the {\it transport
laws}\footnote{Recall from \cite{Br-Spe} that the term {\it
connection} was chosen because of an analogy with path-connections
in differential geometry. In particular, the transport law is a
variation or special case of the transport law for a
path-connection. }.

It is easily verified that the singular cubical complex   $KX$   of  a  space
$X$   satisfies these axioms if  $\circ _j$   is defined by
$$
     (a\circ _j b)(t_1 ,t_2 ,\ldots,t_n ) = \begin{cases}
     a(t_1 ,\ldots, t_{j-1},2t_j,t_{j+1} ,\ldots,t_n) &(t_j \le
     \half)\\
      b(t_1 ,\ldots, t_{j-1},2t_j-1,t_{j+1} ,\ldots,t_n) &(t_j \ge
     \half)\\
     \end{cases}
 $$
whenever   $\del^-_j  b =  \del^+_j  a$.  In this context the
transport law for $\Gamma^-_1 (a\circ b)$ can be illustrated by
the picture \vspace{1in}

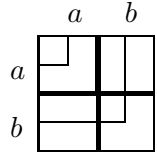
\begin{center}
\setlength{\unitlength}{0.15in}
\begin{picture}(0,0)(9,-1)
 \put(0,0){\framebox(4,4){}} \put(0,1){\line(2,0){3}}
\put(0,2){\line(3,0){2}} \put(2,4){\line(0,-2){2}}
\put(0,3){\line(3,0){1}} \put(2,4){\line(0,-2){2}}
\put(1,4){\line(0,-2){1}} \put(3,4){\line(0,-2){3 }}

\put(0,2){\rule[-0.01in]{0.6in}{0.02in}}
\put(2,0){\rule[0in]{0.02in}{0.6in}}
 \put(1,4.5){\makebox{$a$}}
 \put(3,4.5){\makebox{$b$}}
 \put(-1,0.5){\makebox{$b$}}
\put(-1,2.5){\makebox{$a$}}

\end{picture}
\end{center}

\begin{definition} \label{2.1}{\em A  {\it cubical $\omega$-category
with connections} $G = \{G_n\}$  is a cubical set with connections
and compositions such that each  $\circ _j$   is a category
structure on   $G_n$    with identity elements $\epsilon_j  y \;(y
\in G_{n-1}   )$, and in addition
\begin{equation} \Gamma^+_ix \circ_i\Gamma^-_ix = \eps_{i+1}x, \quad
\Gamma^+_ix \circ_{i+1}\Gamma^-_ix = \eps_{i}x.\tag{2.7}
\end{equation}

For simplicity, a cubical $\omega$-category with connections will
be called a cubical $\omega$-category in the rest of this paper.
}\end{definition}
\begin{remark}{\em  This list is a part of the list of structure and
axioms which first appears in the thesis of Mosa \cite[Chapter
V]{Mo},  in the context of cubical algebroids with connection, and
appears again in the thesis of Al-Agl  \cite{AA}. The rules for
the connections are fairly clear extensions of the axioms given in
\cite{BH1,BH2}, given the general notion of thin structure on a
double category discussed by Spencer in \cite{Spe}. }
\end{remark}

Note that a cubical $\omega$-category has an underlying  cubical
set under its face and degeneracy operations.

It is now straightforward to construct a functor from
$\omega$-categories to cubical $\omega$-categories. The following
type of construction is well known.

\begin{definition} \label{2.2} {\em The {\it cubical nerve\/} of an
$\omega$-category~$X$ is the cubical $\omega$-category  $\lambda
X$ defined as follows: $$(\lambda X)_n=\Hom[M(I^n),X],$$ and the
operations $\d^\alpha_i$, $\epsilon_i$, $\circ _i$,
$\Gamma^{\alpha}_i$ are induced by $\check\d^\alpha_i$,
$\check\epsilon_i$, $\check\mu_i$, $\check\Gamma^{\alpha}_i$
according to the formulae $$\d^\alpha_i
x=x\circ\check\d^\alpha_i\colon M(I^{n-1})\to X$$ for $x\colon
M(I^n)\to X$, etc. }
\end{definition}

In particular, in Definition~\ref{2.2}, note that the domain
of~$\circ _i$ in $(\lambda X)_n\times(\lambda X)_n$ is precisely
$$\Hom[M(I^{i-1}\times[0,2]\times I^{n-i}),X]$$ according to
Remark~\ref{1.10}. To check that $\lambda X$ satisfies the
conditions of Definition~\ref{2.1}, one must check the
corresponding identities for the~$\check\d^\alpha_i$, etc. Many
relations essentially come from properties
of the underlying morphisms $\check\d^\alpha$, etc. The relation
$\d^\alpha_i\epsilon_i=\id$, for example, comes from the easily
checked relation $\check\epsilon\circ\check\d^\alpha=\id$. For
relations involving composition, one must use the morphisms
$\check\iota^\alpha\colon M(I)\to M([0,2])$ which present
$M([0,2])$ as a push-out. Thus, to check the relation
$\d^-_i(x\circ _i y)=\d^-_i x$, which is a relation between binary
operators, one must check that
$$(\check\iota^-)^{-1}\check\mu\check\d^-(\sigma)=\check\d^-(\sigma)$$
and $$(\check\iota^+)^{-1}\check\mu\check\d^-(\sigma)=\emptyset$$
for every cell~$\sigma$ in~$I^0$. For the associative law, one
must consider morphisms from $M(I)$ to $M([0,3])$.

The functoriality of the tensor product is responsible for
formulae looking like commutation rules, such as
$\d^\A_i\eps_j=\eps_{j-1}\d^\A_i$ for $i<j$.

\begin{remark}{\em
Any natural operation $\theta$ on cubical $\omega$-categories
determines an underlying homomorphism $\check{\theta}$ between
$\omega$-categories. For example, if $\theta$ maps $G_n$ to
$G_m$, then in particular $\theta$ maps $$[\lambda
M(I^n)]_n=\Hom[M(I^n),M(I^n)] $$ to $[\lambda
M(I^n)]_m=\Hom[M(I^m),M(I^n)] $ and $\check{\theta}=\theta(\id):
M(I^m) \to M(I^n).$ }\end{remark}

\section{The $\omega$-category associated to a cubical
$\omega$-category}

In this section we construct a functor~$\gamma$ associating an
$\omega$-category to a cubical $\omega$-category. The idea is to
recover an $\omega$-category from its nerve. We will use certain
folding operations, which are defined as follows.

\begin{definition} \label{3.1}{\em  Let $G$ be a cubical $\omega$-category.
The {\it folding operations\/} are the operations
$$\psi_i,\Psi_r,\Phi_m\colon G_n\to G_n$$ defined for $1\leq i\leq
n-1$, $1\leq r\leq n$ and $0\leq m\leq n$ by
\begin{gather*}
 \psi_i x
 =\Gamma^+_i\d^-_{i+1}x\circ _{i+1}x\circ _{i+1}\Gamma^-_i\d^+_{i+1}x,\\
  \Psi_r=\psi_{r-1}\psi_{r-2}\ldots\psi_1,\\  \Phi_m=\Psi_1\Psi_2\ldots\Psi_m
 =\psi_1(\psi_2\psi_1)\ldots(\psi_{m-1}\ldots\psi_1).\end{gather*}
}\end{definition}

Note in particular that $\Psi_1$, $\Phi_0$ and~$\Phi_1$ are
identity operations.

Here is a picture of $\psi_1: G_2 \to G_1$:

\vspace{1.5in}

\begin{center}
\setlength{\unitlength}{0.1in}
\begin{picture}(0,0)(10,-2)
 \put(0,0){\framebox(4,12){$x$}} \put(0,4){\line(1,0){4}}
\put(0,8){\line(1,0){4}}

\put(0,1){\line(1,0){3}}

 \put(0,2){\line(3,0){2}} \put(2,4){\line(0,-2){2}}
\put(0,3){\line(3,0){1}} \put(2,4){\line(0,-2){2}}
\put(1,4){\line(0,-2){1}} \put(3,4){\line(0,-2){3 }}

 \put(1,11){\line(2,0){3}}
\put(1,11){\line(0,-1){3}} \put(1,11){\line(0,-2){2}}
\put(2,10){\line(3,0){2}} \put(2,10){\line(0,-2){2}}
\put(3,9){\line(0,-2){1}} \put(3,9){\line(1,0){1}}
 \put(6,6){\directs{1}{2}}
\put(-7,6){\makebox{$\psi_1(x)= $}}

\end{picture}
\end{center}

The idea behind Definition~\ref{3.1} is best seen from the action
of the underlying endomorphism~$\check\Phi_n$ in the
$\omega$-category of sets $M(I^n)$.

\begin{proposition}\label{3.2} The endomorphism  $\check\Phi_n\colon M(I^n)\to
M(I^n)$ underlying the folding
operation $\Phi_n$ is given by $\check\Phi_n(I^n)=I^n$ and
$$\check\Phi_n(\sigma\times d^\alpha_0 I\times I^p)
 =d^\alpha_p I^n$$
for any cell~$\sigma$ in~$I^{n-p-1}$.
\end{proposition}

\begin{proof}Let $\check\psi\colon M(I^2)\to M(I^2)$ be the
operation underlying~$\psi_1$ in dimension~$2$. The operations
underlying $\psi_i$, $\Psi_r$ and~$\Phi_m$ in dimension~$n$ are
then given by
\begin{gather*}
 \check\psi_i=\id^{i-1}\otimes\check\psi\otimes\id^{n-i-1},\\
  \check\Psi_r=\check\psi_1\check\psi_2\ldots\check\psi_{r-1},\\
   \check\Phi_m=\check\Psi_m\check\Psi_{m-1}\ldots\check\Psi_1.\end{gather*}
One finds that $\check\psi(I^2)=I^2$, from which it follows that
$\check\psi_i(I^n)=I^n$ and then $\check\Phi_n(I^n)=I^n$. One also
finds that
$$\check\psi(\tau\times d^\alpha_0 I)
 =d^\alpha_0 I\times d^\alpha_0 I$$
for any cell~$\tau$ in~$I$. For a cell~$\sigma$ in~$I^{n-p-1}$ it
follows that
$$(\check\Psi_{n-p-1}\ldots\check\Psi_1)
 (\sigma\times d^\alpha_0 I\times I^p)
 \subseteq I^{n-p-1}\times d^\alpha_0 I\times I^p$$
and
$$\check\Psi_{n-p}(\check\Psi_{n-p-1}\ldots\check\Psi_1)
 (\sigma\times d^\alpha_0 I\times I^p)
 =(d^\alpha_0 I)^{n-p}\times I^p.$$
It then follows that $\check\Phi_n(\sigma\times d^\alpha_0 I\times
I^p)$ is independent of~$\sigma$. It now suffices to show that
$$\check\Phi_n[(d^\alpha_0 I)^{n-p}\times I^p]=d^\alpha_p I^n.$$

Recall that $d^\alpha_{n-1}I^n$ is the union of the $(n-1)$-cells
$$\tau_1=d^\alpha_0 I\times I^{n-1},\
 \tau_2=I\times d^{-\alpha}_0 I\times I^{n-2},\ \ldots.$$
We see that
$$\check\Phi_n(\tau_2)
 =\check\Phi_n(d^\alpha_0 I\times d^{-\alpha}_0 I\times I^{n-2})
 \subset\check\Phi_n(\tau_1),$$
etc., so that $\check\Phi(d^\alpha_{n-1}I^n)=\check\Phi(\tau_1)$.
It follows that
$$\check\Phi_n(d^\alpha_0 I\times I^{n-1})=\check\Phi_n(\tau_1)
 =\check\Phi_n(d^\alpha_{n-1}I^n)=d^\alpha_{n-1}\check\Phi_n(I^n)
 =d^\alpha_{n-1}I^n.$$
By similar reasoning,
$$\check\Phi_n[(d^\alpha_0 I)^2\times I^{n-2}]
 =d^\alpha_{n-2}\check\Phi_n(d^\alpha_0 I\times I^{n-1})
 =d^\alpha_{n-2}d^\alpha_{n-1}I^n=d^\alpha_{n-2}I^n,$$
and so on, eventually giving
$$\check\Phi_n[(d^\alpha_0 I)^{n-p}\times I^p]=d^\alpha_p I^n$$
as required. This completes the proof.
\end{proof}

It follows from Proposition~\ref{3.2} that $\check\Phi_n\colon
M(I^n)\to M(I^n)$ is an idempotent endomorphism with image
 $$F_n
 =\{I^n,d^-_{n-1}I^n,d^+_{n-1}I^n,\ldots,d^-_0 I^n,d^+_0 I^n\}.$$
 In fact $F_n$ is nothing else but the $n$-globe.
For an $\omega$-category~$X$, it follows that
 $$\Phi_n[(\lambda X)_n]\cong\Hom(F_n,X).$$
Now, it is clear that $F_n$ has a presentation with
generator~$I^n$ and relations $d^-_n I^n=d^+_n I^n=I^n$; therefore
 $$\Phi_n(\lambda X)_n\cong\{\,x\in X:d^-_n x=d^+_n x=x\}.$$
It follows that $X$ can be recovered from $\lambda X$ as the
colimit of a sequence
 $$\Phi_0[(\lambda X)_0]\to\Phi_1[(\lambda X)_1]\to\cdots.$$
We will now explain how to perform this construction for
cubical $\omega$-categories in general. We begin with some
elementary relations.

\begin{proposition}\label{3.3} The folding operations satisfy the
following relations\/{\rm:}

{\rm (i)}
\begin{alignat*}{2}
 \psi_j\epsilon_i&=\epsilon_i\psi_{j-1}\quad &&{\it for}\ i<j,\\
   \psi_j\epsilon_j&=\psi_j\epsilon_{j+1}=\psi_j\Gamma^{-\alpha}_j
 =\epsilon_j,\quad &&\\
  \psi_j\epsilon_i&=\epsilon_i\psi_j &&{\it for}\ i>j+1,\\
   \d^\alpha_i\psi_j&=\psi_{j-1}\d^\alpha_i &&{\it for}\ i<j,\\
    \d^-_j\psi_j x&=\d^-_j x\circ _j\d^+_{j+1}x,&&\\
     \d^+_j\psi_j x&=\d^-_{j+1}x\circ _j\d^+_j x,&&\\
      \d^\alpha_{j+1}\psi_j&=\epsilon_j\d^\alpha_j\d^\alpha_{j+1},&&\\
      \d^\alpha_i\psi_j&=\psi_j\d^\alpha_i && {\it for}\
      i>j+1;\\
\intertext{\rm (ii)}
              \Psi_1\epsilon_1&=\epsilon_1&&\\
  \Psi_r\epsilon_1&=\epsilon_1\Psi_{r-1}&&{\it for}\ r>1,\\
   \Psi_r\epsilon_i&=\epsilon_{i-1}\Psi_r &&
 {\it for}\ 1<i\leq r,\\
 \d^\alpha_i\Psi_r&=\Psi_r\d^\alpha_i && {\it for}\ i>r,\\
 \d^\alpha_r\Psi_r&=\epsilon_1^{r-1}(\d^\alpha_1)^r;&&\\
\intertext{\rm (iii)}
 \Phi_m\epsilon_i&=\epsilon_1\Phi_{m-1} &&
 {\it for}\ 1\leq i\leq m,\\
 \d^\alpha_i\Phi_m&=\Phi_m\d^\alpha_i &&{\it for}\ i>m,\\
\d^\alpha_m\Phi_m&=\epsilon_1^{m-1}(\d^\alpha_1)^m.&&
\end{alignat*}
\end{proposition}

\begin{proof}(i) These relations are straightforward consequences
of the definitions.

(ii) Since $\Psi_1=\id$, we have $\Psi_1\epsilon_1=\epsilon_1$.

From part~(i), if $r>1$ then
 $$\Psi_r\epsilon_1
 =(\psi_{r-1}\ldots\psi_2)\psi_1\epsilon_1
 =(\psi_{r-1}\ldots\psi_2)\epsilon_1
 =\epsilon_1(\psi_{r-2}\ldots\psi_1)
 =\epsilon_1\Psi_{r-1}.$$

Also from part~(i), if $1<i\leq r$ then
\begin{align*}
 \Psi_r\epsilon_i
 &=(\psi_{r-1}\ldots\psi_i)\psi_{i-1}(\psi_{i-2}\ldots\psi_1)\epsilon_i\\
  &=(\psi_{r-1}\ldots\psi_i)\psi_{i-1}\epsilon_i(\psi_{i-2}\ldots\psi_1)\\
   &=(\psi_{r-1}\ldots\psi_i)\epsilon_{i-1}(\psi_{i-2}\ldots\psi_1)\\
   &=\epsilon_{i-1}(\psi_{r-2}\ldots\psi_{i-1})(\psi_{i-2}\ldots\psi_1)\\
   &=\epsilon_{i-1}\Psi_{r-1}.\end{align*}

From part~(i), if $i>r$ then
$$\d^\alpha_i\Psi_r
 =\d^\alpha_i(\psi_{r-1}\ldots\psi_1)
 =(\psi_{r-1}\ldots\psi_1)\d^\alpha_i
 =\Psi_r\d^\alpha_i.$$

It now follows that
\begin{align*}
 \d^\alpha_r\Psi_r
 &=\d^\alpha_r\psi_{r-1}\Psi_{r-1}\\
   &=\epsilon_{r-1}\d^\alpha_{r-1}\d^\alpha_r\Psi_{r-1}\\
   &=\epsilon_{r-1}\d^\alpha_{r-1}\Psi_{r-1}\d^\alpha_r\\
    &=\ldots\\  &=\epsilon_{r-1}\ldots\epsilon_2\epsilon_1\d^\alpha_1\Psi_1
 \d^\alpha_2\ldots\d^\alpha_r\\
  &=\epsilon_{r-1}\ldots\epsilon_2\epsilon_1\d^\alpha_1
 \d^\alpha_2\ldots\d^\alpha_r\\
  &=\epsilon_1^{r-1}(\d^\alpha_1)^r,\end{align*}
using (2.1).

(iii) From part~(ii), if $1\leq i\leq m$ then
\begin{align*}
 \Phi_m\epsilon_i
 &=\Psi_1(\Psi_2\ldots\Psi_{m-i+1})(\Psi_{m-i+2}\ldots\Psi_m)\epsilon_i\\
  &=\Psi_1(\Psi_2\ldots\Psi_{m-i+1})\epsilon_1
 (\Psi_{m-i+1}\ldots\Psi_{m-1})\\
  &=\Psi_1\epsilon_1(\Psi_1\ldots\Psi_{m-i})(\Psi_{m-i+1}\ldots\Psi_{m-1})\\
   &=\epsilon_1(\Psi_1\ldots\Psi_{m-i})(\Psi_{m-i+1}\ldots\Psi_{m-1})\\
    &=\epsilon_1\Phi_{m-1}.\end{align*}

Also from part~(ii), if $i>m$ then
$$\d^\alpha_i\Phi_m
 =\d^\alpha_i(\Psi_1\ldots\Psi_m)
 =(\Psi_1\ldots\Psi_m)\d^\alpha_i
 =\Phi_m\d^\alpha_i.$$

It now follows that \begin{align*}\d^\alpha_m\Phi_m
 &=\d^\alpha_m\Phi_{m-1}\Psi_m\\
 &=\Phi_{m-1}\d^\alpha_m\Psi_m \\
 &=\Phi_{m-1}\epsilon_1^{m-1}(\d^\alpha_1)^m\\
 &=\epsilon_1^{m-1}\Phi_0(\d^\alpha_1)^m\\
 &=\epsilon_1^{m-1}(\d^\alpha_1)^m.
 \end{align*}
\end{proof}

We now observe that the operators~$\psi_i$ are idempotent, and
characterise their images.

\begin{proposition}\label{3.4} Let $G$ be a cubical $\omega$-category, and
suppose that $1\leq i\leq n-1$. The operator  $\psi_i\colon
G_n\to G_n$ is idempotent. An element~$x$ of~$G_n$ is in
$\psi_i(G_n)$ if and only if $\d^-_{i+1}x$ and $\d^+_{i+1}x$ are
in $\Im\epsilon_i$.
\end{proposition}

\begin{proof}From Proposition 3.3(i), if $x\in\psi_i(G_n)$ then
$\d^-_{i+1}x$ and $\d^+_{i+1}x$ are in $\Im\epsilon_i$.

To complete the proof, suppose that $\d^-_{i+1}x$ and
$\d^+_{i+1}x$ are in $\Im\epsilon_i$; it suffices to show that
$\psi_i x=x$. Now,
$$\Gamma^{-\A}_i\d^\alpha_{i+1}x\in\Im\Gamma^{-\alpha}_i\epsilon_i
 =\Im\epsilon_i^2=\Im\epsilon_{i+1}\epsilon_i,$$
so that the $\Gamma^{-\alpha}_i\d^\alpha_{i+1}x$ are identities
for~$\circ _{i+1}$. It follows that
$$\psi_i x
 =\Gamma^+_i\d^-_{i+1}x\circ _{i+1}x\circ _{i+1}\Gamma^-_i\d^+_{i+1}x
 =x$$
as required. This completes the proof.
\end{proof}

There is a similar result for~$\Phi_n$ as follows.

\begin{proposition} \label{3.5} Let $G$ be a cubical $\omega$-category. The
operator \newline $\Phi_n\colon G_n\to G_n$ is idempotent. An
element $x$ of~$G_n$ is in $\Phi_n(G_n)$ if and only if
$\d^\alpha_m x\in\Im\epsilon_1^{m-1}$ for $1\leq m\leq n$ and
$\alpha=\pm$.
\end{proposition}

\begin{proof}Since $\Phi_n=\Phi_m(\Psi_{m+1}\ldots\Psi_n)$, it
follows from Proposition 3.3(iii) that
$$\Im\d^\alpha_m\Phi_n\subset\Im\d^\alpha_m\Phi_m
 \subset\Im\epsilon_1^{m-1}.$$
Conversely, suppose that $\d^\alpha_m x\in\Im\epsilon_1^{m-1}$ for
$1\leq m\leq n$ and $\alpha=\pm$; it suffices to show that $\Phi_n
x=x$, and for this it suffices to show that $\psi_i x=x$ for
$1\leq i\leq n-1$. But
$$\d^\alpha_{i+1}x\in\Im\epsilon_1^i=\Im\epsilon_i\epsilon_1^{i-1}
 \subset\Im\epsilon_i$$
for $\alpha=\pm$, so that $\psi_i x=x$ by Proposition~\ref{3.4}. This
completes the proof.
\end{proof}

There is a useful result related to Proposition~\ref{3.5} as
follows.

\begin{proposition} \label{3.6} If $x\in\Phi_n(G_n)$ and $1\leq m\leq
n$, then
$$\d^\alpha_m x=\epsilon_1^{m-1}(\d^\alpha_1)^m x\quad
 {\it and}\quad
 \epsilon_m\d^\alpha_m x=\epsilon_1^m(\d^\alpha_1)^m x.$$
\end{proposition}

\begin{proof}By Proposition~\ref{3.5}, $\d^\alpha_m
x=\epsilon_1^{m-1}x'$ for some~$x'$. It follows that
$$x'
 =(\d^\alpha_1)^{m-1}\epsilon_1^{m-1}x'
 =(\d^\alpha_1)^{m-1}\d^\alpha_m x
 =(\d^\alpha_1)^m x,$$
so that $\d^\alpha_m
x=\epsilon_1^{m-1}x'=\epsilon_1^{m-1}(\d^\alpha_1)^m x$ and
$\epsilon_m\d^\alpha_m x=\epsilon_m\epsilon_1^{m-1}(\d^\alpha_1)^m
x=\epsilon_1^m(\d^\alpha_1)^m x$ as required.
\end{proof}

We now deduce various closure properties for the family of sets
$\Phi_n(G_n)$.

\begin{proposition} \label{3.7} Let $G$ be a cubical $\omega$-category. The
family of sets $\Phi_n(G_n)$
{\rm(}$n\geq 0${\rm)} is closed under the~$\d^\alpha_i$ and
under~$\epsilon_1$. The individual sets $\Phi_n(G_n)$ are closed
under $\epsilon_i\d^\alpha_i$ and~$\circ _i$ for $1\leq i\leq n$.
\end{proposition}

\begin{proof} We use the characterisation in Proposition~\ref{3.5}.
We first show that the family is closed under~$\d^\alpha_1$.
Indeed, if $x\in\Phi_n(G_n)$ then
 $$\d^\beta_m\d^\alpha_1 x
 =\d^\alpha_1\d^\beta_{m+1}x
 \in\d^\alpha_1(\Im\epsilon_1^m)
 =\Im\epsilon_1^{m-1},$$
since $\d^\alpha_1\epsilon_1=\id$.

Next we show that the family is closed under~$\epsilon_1$. Indeed,
if $x\in\Phi_n(G_n)$, then $\d^\beta_1\epsilon_1
x\in\Im\epsilon_1^0$ trivially, and for $m>1$ we have
 $$\d^\beta_m\epsilon_1 x
 =\epsilon_1\d^\beta_{m-1}x
 \in\epsilon_1(\Im\epsilon_1^{m-2})
 =\Im\epsilon_1^{m-1}.$$

It now follows from Proposition~\ref{3.6} that the family is
closed under~$\d^\alpha_i$ for all~$i$. Similarly, $\Phi_n(G_n)$
is closed under $\epsilon_i\d^\alpha_i$.

It remains to show that $x\circ _i y\in\Phi_n(G_n)$ when
$x$~and~$y$ are in $\Phi_n(G_n)$ and the composite exists. Suppose
that $\d^\beta_m x=\epsilon_1^{m-1}x'$ and $\d^\beta_m
y=\epsilon_1^{m-1}y'$. If $m<i$ then
 \begin{gather*}\d^\beta_m(x\circ _i y)
 =\d^\beta_m x\circ _{i-1}\d^\beta_m y
 =\epsilon_1^{m-1}x'\circ _{i-1}\epsilon_1^{m-1}y'\\
 =\epsilon_1^{m-1}(x'\circ _{i-m}y')
 \in\Im\epsilon_1^{m-1};\end{gather*}
if $m=i$ then $\d^\beta_m(x\circ _i y)$ is $\d^\beta_m
x=\epsilon_1^{m-1}x'$ or $\d^\beta_m y=\epsilon_1^{m-1}y'$, so
$\d^\beta_m(x\circ _i y)$ is certainly in $\Im\epsilon_1^{m-1}$;
and if $m>i$ then
 \begin{gather*}\d^\beta_m(x\circ _i y)
 =\d^\beta_m x\circ _i\d^\beta_m y
 =\epsilon_1^{m-1}x'\circ _i\epsilon_1^{m-1}y'  =\\
 =\epsilon_1^{i-1}(\epsilon_1^{m-i}x'\circ _1\epsilon_1^{m-i}y')
 =\epsilon_1^{i-1}\epsilon_1^{m-i}x'
 \in\Im\epsilon_1^{m-1}\end{gather*}
(note that $\epsilon_1^{m-i}y'$ is an identity for~$\circ _1$
because it lies in the image of~$\epsilon_1$).

This completes the proof.
\end{proof}

We now obtain the desired sequence of $\omega$-categories.

\begin{theorem} \label{3.8} Let $G$ be a cubical $\omega$-category. Then
there is a sequence of $\omega$-categories and homomorphisms
$$\xymatrix@1{
 \Phi_0(G_0)\ar[r]^{\epsilon_1} &
 \Phi_1(G_1)\ar[r]^{\epsilon_1} &
 \Phi_2(G_2)\ar[r] &
 \cdots}$$
with the following structure on $\Phi_n(G_n)${\rm:} if $0\leq p<n$
then $$d^\alpha_p x=(\epsilon_1)^{n-p}(\d^\alpha_1)^{n-p}x$$ and
$x\c_p y=x\circ _{n-p}y$ where defined\/{\rm;} if $p\geq n$ then
$d^\alpha_p x=x$ and the only composites are given by $x\c_p x=x$.
\end{theorem}

\begin{proof} We first show that for a fixed value of~$n$ the given
structure maps $d^\alpha_p$ and~$\c_p$ make $\Phi_n(G_n)$ into an
$\omega$-category. By Proposition~\ref{3.7}, $\Phi_n(G_n)$ is
closed under the structure maps for $0\leq p<n$, and the same
result holds trivially for $p\geq n$.

From the identities in Section 2, if $0\leq p<n$ then the triple
$$(d^-_p,d^+_p,\c_p)
 =(\epsilon_{n-p}\d^-_{n-p},\epsilon_{n-p}\d^+_{n-p},\circ _{n-p})$$
makes $\Phi_n(G_n)$ into the morphism set of a category (with
$d^-_p x$ and $d^+_p x$ the left and right identities of~$x$ and
with $\c_p$ as composition), and  these structures commute with
one another. Trivially the triples $(d^-_p,d^+_p,\c_p)$ for
$n\geq p$ provide further commuting category structures. To show
that these structures make $\Phi_n(G_n)$ into an
$\omega$-category, it now suffices to show that an identity
for~$\c_p$ is also an identity for~$\c_q$ if $q>p$; in other
words, it suffices to show that $d^\beta_q d^\alpha_p
x=d^\alpha_p x$ for $x\in\Phi_n(G_n)$ and $q>p$. For $q\geq n$,
this is trivial; we may therefore assume that $0\leq p<q<n$. But
Proposition~\ref{3.6} gives us \begin{align*}d^\beta_q d^\alpha_p
x
 &=\epsilon_1^{n-q}(\d^\beta_1)^{n-q}
 \epsilon_1^{n-p}(\d^\alpha_1)^{n-p}x\\
 &=\epsilon_1^{n-q}\epsilon_1^{q-p}(\d^\alpha_1)^{n-p}x\\
 &=\epsilon_1^{n-p}(\d^\alpha_1)^{n-p}x\\
 &=d^\alpha_p x\end{align*}
as required.

We have now shown that the $\Phi_n(G_n)$ are $\omega$-categories.
We know from Proposition~\ref{3.7} that $\epsilon_1$ maps
$\Phi_n(G_n)$ into $\Phi_{n+1}(G_{n+1})$, and it remains to show
that this function is a homomorphism. That is to say, we must show
that $\epsilon_1 d^\alpha_p x=d^\alpha_p\epsilon_1 x$ for
$x\in\Phi_n(G_n)$, and we must show that $\epsilon_1(x\c_p
y)=\epsilon_1 x\c_p\epsilon_1 y$ for $x\c_p y$ a composite in
$\Phi_n(G_n)$. But if $0\leq p<n$ then $$\epsilon_1 d^\alpha_p x
 =\epsilon_1\epsilon_{n-p}\d^\alpha_{n-p}x
 =\epsilon_{n-p+1}\d^\alpha_{n-p+1}\epsilon_1 x
 =d^\alpha_p\epsilon_1 x$$
and
$$\epsilon_1(x\c_p y)
 =\epsilon_1(x\circ _{n-p}y)
 =\epsilon_1 x\circ _{n-p+1}\epsilon_1 y
 =\epsilon_1 x\c_p\epsilon_1 y$$
by identities in Section 2; if $p=n$ we get
$$\epsilon_1 d^\alpha_n x
 =\epsilon_1 x
 =\epsilon_1\d^\alpha_1\epsilon_1 x
 =d^\alpha_n\epsilon_1 x$$
and
$$\epsilon_1(x\c_n x)
 =\epsilon_1 x
 =\epsilon_1 x\circ _1\epsilon_1 x
 =\epsilon_1 x\c_n\epsilon_1 x;$$
and if $p>n$ then
$$\epsilon_1 d^\alpha_p x=\epsilon_1 x=d^\alpha_p\epsilon_1 x$$
and
$$\epsilon_1(x\c_p x)=\epsilon_1 x=\epsilon_1 x\c_p\epsilon_1 x$$
trivially. This completes the proof.
\end{proof}

We can now define a functor from cubical $\omega$-categories to
$\omega$-categories.

\begin{definition} \label{3.9}{\em  Let $G$ be a cubical $\omega$-category.
The $\omega$-category $\gamma G$
{\it associated to~$G$} is the colimit of the sequence
$$\xymatrix@1{
 \Phi_0(G_0)\ar[r]^{\epsilon_1} &
 \Phi_1(G_1)\ar[r]^{\epsilon_1} &
 \Phi_2(G_2)\ar[r] &
 \cdots.}$$}
\end{definition}

\begin{remark} \label{3.10}{\em  In Definition~\ref{3.9}, one can identify
$\Phi_n(G_n)$ with the subset of $\gamma G$ consisting of
elements~$x$ such that $d^-_n x=d^+_n x=x$. Indeed,
the~$\epsilon_1$ are injective, because
$\d^\alpha_1\epsilon_1=\id$, so that $\Phi_n(G_n)$ can be
identified with a subset of $\gamma G$; if $x\in\Phi_n(G_n)$ then
$d^-_n x=d^+_n x=x$ by Theorem~\ref{3.8}; if $x\in\Phi_m(G_m)$ with
$m>n$ and $d^-_n x=d^+_n x=x$, then
$$x=\epsilon_{m-n}\d^-_{m-n}x=\epsilon_1^{m-n}(\d^-_1)^{m-n}x$$
(Proposition~\ref{3.6}) with $(\d^-_1)^{m-n}x\in\Phi_n(G_n)$
(Proposition~\ref{3.7}), and $x$ can be identified with
$(\d^-_1)^{m-n}x$.}
\end{remark}
\begin{remark}{\em  It is convenient to describe the $\omega$-category $\gamma
G$ in terms of the folding operations, but one can get a more
direct description by using Proposition 3.5. The more direct
description needs face maps, degeneracies and compositions, but
not connections.}
\end{remark}
\section{The natural isomorphism $A\colon\gamma\lambda X\to
X$}

Let $X$ be an $\omega$-category. From Definition~\ref{2.2} there
is a cubical  $\omega$-category $\lambda X$, and from
Definition~\ref{3.9} there is an $\omega$-category $\gamma\lambda
X$. We will now construct a natural isomorphism
$A\colon\gamma\lambda X\to X$.

Let $F_n$ be the $\omega$-category with one generator~$I^n$ and
with relations $d^-_n I^n=d^+_n I^n=I^n$. By Proposition~\ref{3.2},
$F_n$~can be realised as a sub-$\omega$-category of $M(I^n)$, and
the morphism $\check\Phi_n\colon M(I^n)\to M(I^n)$ associated to
the folding operation~$\Phi_n$ is an idempotent operation with
image equal to~$F_n$. Recalling that $(\lambda
X)_n=\Hom[M(I^n),X]$, we see that
$$\Phi_n[(\lambda X)_n]
 =\{\,x\in\Hom[M(I^n),X]:x\check\Phi_n=x\,\}.$$
Let
$$A\colon\Phi_n[(\lambda X)_n]\to X$$
be the function given by
$$A(x)=x(I^n);$$
we see that $A$ is an injection with image equal to
$$\{\,x\in X:d^-_n x=d^+_n x=x\,\}.$$
These functions are compatible with the sequence
$$\xymatrix@1{
 \cdots\ar[r] &
 \Phi_n[(\lambda X)_n]\ar[r]^-{\epsilon_1} &
 \Phi_{n+1}[(\lambda X)_{n+1}]\ar[r] &
 \cdots;}$$
indeed, if $x\in\Phi_n[(\lambda X)_n]$ then
$$A(\epsilon_1 x)
 =(\epsilon_1 x)(I^{n+1})
 =x\check\epsilon_1(I^{n+1})
 =x(I^n)
 =A(x).$$
The functions $A\colon\Phi_n[(\lambda X)_n]\to X$ therefore induce
a bijection $A\colon\gamma\lambda X\to X$. We will now prove the
following result.

\begin{theorem} \label{4.1} The functions $A\colon\gamma\lambda X\to X$
form a natural isomorphism of $\omega$-categories.
\end{theorem}

\begin{proof}We have already shown that $A\colon\gamma\lambda X\to
X$ is a bijection, and it is clearly natural. It remains to show
that $A$~is a homomorphism. It suffices to show that
$$A\colon\Phi_n[(\lambda X)_n]\to X$$
is a homomorphism for each~$n$; in other words, we must show that
$A(d^\alpha_p x)=d^\alpha_p A(x)$ for $x\in \Phi_n[(\lambda X)_n]$
and that $A(x\c_p y)=A(x)\c_p A(y)$ for $x\c_p y$ a composite in
$\Phi_n[(\lambda X)_n]$.

Suppose that $x\in\Phi_n[(\lambda X)_n]$ and $0\leq p<n$. Noting
that $x=x\check\Phi_n$ and using Proposition~\ref{3.2}, we find that
\begin{align*}
 A(d^\alpha_p x)
 &=A(\epsilon_{n-p}\d^\alpha_{n-p}x)\\
 &=(\epsilon_{n-p}\d^\alpha_{n-p}x)(I^n)\\
 &=x\check\d^\alpha_{n-p}\check\epsilon_{n-p}(I^n)\\
 &=x(I^{n-p-1}\times d^\alpha_0 I\times I^p)\\
 &= x\check\Phi_n(I^{n-p-1}\times d^\alpha_0 I\times I^p)\\
 &=x(d^\alpha_p I^n)\\  &=d^\alpha_p x(I^n)\\
  &=d^\alpha_p A(x).\end{align*}

Suppose that  $x\in\Phi_n[(\lambda X)_n]$ and $p\geq n$. Then
$$A(d^\alpha_p x)
 =A(x)
 =x(I^n)
 =x(d^\alpha_p I^n)
 =d^\alpha_p x(I^n)
 =d^\alpha_p x.$$

Suppose that $x\# _p y$ is a composite in $\Phi_n[(\lambda X)_n]$
with $0\leq p<n$. Let
$$(x,y)\colon  M(I^{n-p-1}\times[0,2]\times I^p)\to X$$
be the morphism such that $(x,y)\check\iota^-_{n-p}=x$ and
$(x,y)\check\iota^+_{n-p}=y$; then
$$A(x\c_p y)=A(x\circ _{n-p}y)=(x,y)\check\mu_{n-p}(I^n).$$
Let $\eta\colon F_n\to M(I^n)$ be the inclusion and let $\pi\colon
M(I^n)\to F_n$ be $\check\Phi_n$ with its codomain restricted
to~$F_n$, so that $\check\Phi_n=\eta\pi$. Since $x$ and~$y$ are in
$\Phi_n[(\lambda X)_n]$, we have $x\check\Phi_n=x$ and
$y=y\check\Phi_n$; we therefore get
$$(x,y)\check\mu_{n-p}(I^n)=(x\eta\pi,y\eta\pi)\check\mu_{n-p}(I^n).$$
Now let $F_p$ be the $\omega$-category with one generator~$z$ and
with relations $d^-_p z=d^+_p z=z$. We see that there is a
factorisation
$$(x\eta\pi,y\eta\pi)=(x\eta,y\eta)(\pi,\pi)$$
through the obvious push-out of
$$F_n\leftarrow F_p\rightarrow F_n.$$
We also see that
$$(\pi,\pi)\check\mu_{n-p}(I^n)
 =\pi\check\iota^-_{n-p}(I^n)\c_p\pi\check\iota^+_{n-p}(I^n).$$
It now follows that
\begin{align*}
 (x\eta\pi,y\eta\pi)\check\mu_{n-p}(I^n)
 &=(x\eta,y\eta)(\pi,\pi)\check\mu_{n-p}(I^n)\\
  &=(x\eta,y\eta)
 [\pi\check\iota^-_{n-p}(I^n)\c_p\pi\check\iota^+_{n-p}(I^n)]\\
   &=(x\eta,y\eta)\pi\check\iota^-_{n-p}(I^n)
 \c_p(x\eta,y\eta)\pi\check\iota^+_{n-p}(I^n)\\
  &=x\eta\pi(I^n)\c_p y\eta\pi(I^n)\\
   &=x(I^n)\c_p y(I^n)\\  &=A(x)\c_p A(y);\end{align*}
therefore
$$A(x\c_p y)=A(x)\c_p A(y).$$

Finally, suppose that $x\# _p y$ is a composite in
$\Phi_n[(\lambda X)_n]$ with $p\geq n$. We must have $x=y$, and we
get
$$A(x\c_p x)=A(x)=A(x)\c_p A(x).$$

This completes the proof.
\end{proof}

\section{Foldings, degeneracies and connections}

According to Theorem~\ref{4.1}, there are natural isomorphisms
$A\colon \gamma\lambda X\to X$ for $\omega$-categories $X$. To
prove that $\omega$-categories are equivalent to cubical
$\omega$-categories, we will eventually construct natural
isomorphisms $B\colon G\to\lambda\gamma G$ for cubical
$\omega$-categories $G$. We will need properties of the folding
operations, and we now begin to describe these.

We first show that the operations~$\psi_i$ behave like the
standard generating transpositions of the symmetric groups (except
of course that they are idempotent rather than involutory, by
Proposition~\ref{3.4}). There are two types of relation, the first of
which is easy.

\begin{proposition}\label{5.1} If $|i-j|\geq 2$, then
$$\psi_i\psi_j=\psi_j\psi_i .$$
\end{proposition}

\begin{proof}This follows from the identities in Section 2.
\end{proof}

The next result is harder.

\begin{theorem} \label{5.2} If $i>1$ then
$$\psi_i\psi_{i-1}\psi_i=\psi_{i-1}\psi_i\psi_{i-1}.$$
\end{theorem}

\begin{proof}Recall the  matrix notation used for certain composites:
if
$$(a_{11}\circ _i\cdots\circ _i a_{1n})
 \circ _{i+1}\cdots\circ _{i+1}
 (a_{m1}\circ _i\cdots\circ _i a_{mn})$$
and
$$(a_{11}\circ _{i+1}\cdots\circ _{i+1} a_{m1})
 \circ _i\cdots\circ _i
 (a_{1n}\circ _{i+1}\cdots\circ _{i+1} a_{mn})$$
are equal by the interchange law, then we will write
$$\begin{bmatrix}
 a_{11}& \ldots& a_{1n}\\  \vdots&&\vdots\\
  a_{m1}& \ldots& a_{mn}  \end{bmatrix} \quad \directs{i}{i+1}$$
for the common value. In such a matrix, we write $-$ for elements
in the image of~$\epsilon_i$ (which are the identities
for~$\circ _i$), and we write $|$ for elements in the image
of~$\epsilon_{i+1}$ (which are the identities for~$\circ _{i+1}$).

We first compute $\psi_i\psi_{i-1}\psi_i x$. It is straightforward
to check that
\begin{align*}\psi_{i-1}\psi_i x&=
 \begin{bmatrix}
 |& \Gamma^+_i\d^-_{i+1}x& \Gamma^-_{i-1}\d^-_{i+1}x\\
 \Gamma^+_{i-1}\d^-_i x& x& \Gamma^-_{i-1}\d^+_i x\\
  \Gamma^+_{i-1}\d^+_{i+1}x& \Gamma^-_i\d^+_{i+1}x& | \end{bmatrix}
  \quad \directs{i}{i+1}.\\
\intertext{It follows that}
 \Gamma^+_i\d^-_{i+1}\psi_{i-1}\psi_i x
 &=\Gamma^+_i(\Gamma^+_{i-1}\d^-_i\d^-_i x
 \circ _i\epsilon_i\d^-_i\d^-_i x
 \circ _i\Gamma^-_{i-1}\d^-_i\d^-_i x)\\
 &=\Gamma^+_i(\Gamma^+_{i-1}\d^-_i\d^-_i x
 \circ _i\Gamma^-_{i-1}\d^-_i\d^-_i x)\\
 &=\Gamma^+_i\epsilon_{i-1}\d^-_i\d^-_i x\\
  &=\epsilon_{i-1}\Gamma^+_{i-1}\d^-_i\d^-_i x\\
\intertext{and}
\Gamma^-_i\d^+_{i+1}\psi_{i-1}\psi_i x
  &=\epsilon_{i-1}\Gamma^-_{i-1}\d^+_i\d^+_i x; \\
\intertext{therefore}
\psi_i\psi_{i-1}\psi_i x
 &=\epsilon_{i-1}\Gamma^+_{i-1}\d^-_i\d^-_i x
 \circ _{i+1}\psi_{i-1}\psi_i x
 \circ _{i+1}\epsilon_{i-1}\Gamma^-_{i-1}\d^+_i\d^+_i x. \end{align*}

Similarly, $\psi_{i-1}\psi_i\psi_{i-1}x$ is as a composite
\begin{small} \directs{i}{i+1} \end{small}
\begin{scriptsize}
$$\begin{bmatrix}
 -& \Gamma^+_{i-1}\Gamma^+_{i-1}\d^-_i\d^-_i x& -& -&
 \Gamma^-_{i-1}\Gamma^+_{i-1}\d^-_i\d^-_i x\\
  -& |& \Gamma^+_i\d^-_{i+1}x& -& \Gamma^-_{i-1}\d^-_{i+1}x\\
   -& |& |& \Gamma^+_i\Gamma^-_{i-1}\d^-_i\d^+_i x&
 \Gamma^-_i\Gamma^-_{i-1}\d^-_i\d^+_i x\\
  -& \Gamma^+_{i-1}\d^-_i x& x& \Gamma^-_{i-1}\d^+_i x& -\\
  \Gamma^+_i\Gamma^+_{i-1}\d^+_i\d^-_i x&
 \Gamma^-_i\Gamma^+_{i-1}\d^+_i\d^-_i x& |& |& -\\
  \Gamma^+_{i-1}\d^+_{i+1}x& -& \Gamma^-_i\d^+_{i+1}x& |& -\\
   \Gamma^+_{i-1}\Gamma^-_{i-1}\d^+_i\d^+_i x& -& -&
 \Gamma^-_{i-1}\Gamma^-_{i-1}\d^+_i\d^+_i x& -\\  \end{bmatrix} $$
\end{scriptsize}
We now evaluate the rows of the matrix for
$\psi_{i-1}\psi_i\psi_{i-1}x$. The first and last rows yield
$\epsilon_{i-1}\Gamma^+_{i-1}\d^-_i\d^-_i x$ and
$\epsilon_{i-1}\Gamma^-_{i-1}\d^+_i\d^+_i x$. The composite of the
non-identity elements in the third row is
$\epsilon_{i+1}\Gamma^-_{i-1}\d^-_i\d^+_i x$, which is an identity
for~$\circ _{i+1}$, so the third row can be omitted. Similarly, the
fifth row can be omitted. The second, fourth and sixth rows have
the same values as the rows of the matrix for $\psi_{i-1}\psi_i
x$. It follows that
$$\psi_{i-1}\psi_i\psi_{i-1}x
 =\epsilon_{i-1}\Gamma^+_{i-1}\d^-_i\d^-_i x
 \circ _{i+1}\psi_{i-1}\psi_i x
 \circ _{i+1}\epsilon_{i-1}\Gamma^-_{i-1}\d^+_i\d^+_i x$$
also. Therefore $\psi_i\psi_{i-1}\psi_i
x=\psi_{i-1}\psi_i\psi_{i-1}x$. This completes the proof.
\end{proof}

\begin{remark}\label{5.3}{ \em Proposition~\ref{5.1} and Theorem~\ref{5.2}
in some sense explain the formula for~$\Phi_m$ in
Definition~\ref{3.1}. The~$\psi_i$ behave like the generating
transpositions $(i,i+1)$ in the symmetric group of permutations of
$\{1,\ldots,m\}$, and, as in the symmetric group, there are
$m!$~distinct composites of $\psi_1,\ldots,\psi_{m-1}$ given by
$$\Psi_{1,l(1)}\ldots\Psi_{m,l(m)}$$ for $1\leq l(r)\leq r$, where
$\Psi_{r,l(r)}=\psi_{r-1}\psi_{r-2}\ldots\psi_{l(r)}$.
 }\end{remark}
Note that $\partial^\pm_1 \Phi_m$ is the hemispherical
decomposition of the $m$-categorical $m$-cube. This orders the
top dimensional cells in the upper hemisphere in reverse order to
the top dimensional cells of the lower hemisphere. Thus the
composite $\Phi_m$ corresponds to the order-reversing permutation
$p \mapsto m+1-p$. In our context, we can characterise $\Phi_m$ as
the zero element in the semigroup generated by $\psi_1,
\ldots,\psi_{m-1}$ as follows.

\begin{theorem} \label{5.4} If $1\leq i\leq m-1$
then
$$\Phi_m\psi_i=\Phi_m.$$
\end{theorem}

\begin{proof}If $r>1$ then
$$\Psi_r\psi_1=(\psi_{r-1}\ldots\psi_2)\psi_1\psi_1
 =(\psi_{r-1}\ldots\psi_2)\psi_1=\Psi_r,$$
since $\psi_1$~is idempotent by Proposition~\ref{3.4}. For $1<i<r$, it
follows from Proposition~\ref{5.1} and Theorem~\ref{5.2} that
\begin{align*}
 \Psi_r\psi_i
 &=(\psi_{r-1}\ldots\psi_{i+1})\psi_i\psi_{i-1}(\psi_{i-2}\ldots\psi_1)
 \psi_i\\  &=(\psi_{r-1}\ldots\psi_{i+1})\psi_i\psi_{i-1}\psi_i
 (\psi_{i-2}\ldots\psi_1)\\  &=(\psi_{r-1}\ldots\psi_{i+1})\psi_{i-1}\psi_i\psi_{i-1}
 (\psi_{i-2}\ldots\psi_1)\\  &=\psi_{i-1}(\psi_{r-1}\ldots\psi_{i+1})\psi_i\psi_{i-1}
 (\psi_{i-2}\ldots\psi_1)\\  &=\psi_{i-1}\Psi_r.\end{align*}
For $1\leq i\leq m-1$, it now follows that
\begin{align*}
 \Phi_m\psi_i
 &=\Psi_1(\Psi_2\ldots\Psi_{m-i})\Psi_{m-i+1}
 (\Psi_{m-i+2}\ldots\Psi_m)\psi_i\\  &=\Psi_1(\Psi_2\ldots\Psi_{m-i})\Psi_{m-i+1}\psi_1
 (\Psi_{m-i+2}\ldots\Psi_m)\\  &=\Psi_1(\Psi_2\ldots\Psi_{m-i})\Psi_{m-i+1}
 (\Psi_{m-i+2}\ldots\Psi_m)\\  &=\Phi_m,\end{align*}
as required.
\end{proof}

We can now give some interactions between $\Phi_m$, degeneracies
and connections. First we have the following result.

\begin{proposition}\label{5.5} For all~$i$ there is a relation
$$\psi_i\Gamma^{\alpha}_i=\epsilon_i.$$
\end{proposition}

\begin{proof}From the definitions we get
\begin{align*}
 \psi_i\Gamma^+_i x
 &=\Gamma^+_i\d^-_{i+1}\Gamma^+_i x\circ _{i+1}\Gamma^+_i x
 \circ _{i+1}\Gamma^-_i\d^+_{i+1}\Gamma^+_i x\\
   &=\Gamma^+_i\epsilon_i\d^-_i x\circ _{i+1}\Gamma^+_i x
 \circ _{i+1}\Gamma^-_i x\\
   &=\epsilon_i^2\d^-_i x\circ _{i+1}\epsilon_i x\\
   &=\epsilon_{i+1}\epsilon_i\d^-_i x\circ _{i+1}\epsilon_i x\\
    &=\epsilon_i x, \end{align*}
and we similarly get $\psi_i\Gamma^-_i x=\epsilon_i x$.
\end{proof}

We draw the following conclusions.

\begin{theorem} \label{5.6} If $1\leq i\leq m$ then
$$\Phi_m\epsilon_i=\epsilon_1\Phi_{m-1}.$$
If $1\leq i\leq m-1$ then
$$\Phi_m\Gamma^{\alpha}_i=\epsilon_1\Phi_{m-1}.$$
\end{theorem}

\begin{proof}The first of these results was given in Proposition
3.3(iii). The second result then follows from Theorem~\ref{5.4} and
Proposition~\ref{5.5}: indeed, we get
$$\Phi_m\Gamma^{\alpha}_i=\Phi_m\psi_i\Gamma^{\alpha}_i
 =\Phi_m\epsilon_i=\epsilon_1\Phi_{m-1}$$
as required.
\end{proof}

\section{Foldings, face maps and compositions}

In this section we describe interactions between the~$\Phi_n$,
face maps and compositions. For face maps, the basic results are
given in Proposition~\ref{3.3}. For compositions, the basic results are
as follows, of which the first two cases correspond to the
2-dimensional case in \cite[Proposition 5.1]{Br-Mo}.

\begin{proposition}\label{6.1} In a cubical $\omega$-category
$$ \psi_i(x\circ _j y) = \begin{cases} (\psi_i x\circ _{i+1}\epsilon_i\d^+_{i+1}y)
 \circ _i(\epsilon_i\d^-_{i+1}x\circ _{i+1}\psi_i y)& \text{if }
 j=i, \\
(\epsilon_i\d^-_i x\circ _{i+1}\psi_i y)
 \circ _i(\psi_i x\circ _{i+1}\epsilon_i\d^+_i y) & \text{if }
 j=i+1,\\
 \psi_i x\circ _j\psi_i y  & \text{otherwise.}
\end{cases}
$$
\end{proposition}
\begin{proof}
Note that we have $\d^+_j x= \d^-_j y $ for $x\circ _j y$ to be
defined.

The proof for the cases $j=i$ and $j=i+1$ consists in evaluating
in two ways each of the matrices $$ \begin{bmatrix} \Gamma^+_i
\d^-_{i+1}x & \eps_i\d^-\dip x \\ \eps \dip \d ^-\dip x & \Gamma
^+_i \d ^-\dip y \\ x & y \\ \Gamma^-_i \d ^+\dip x & \eps\dip
\d^+\dip y \\ \eps_i \d^+\dip y & \Gamma^-_i  \d^+\dip y
\end{bmatrix} \qquad
 \begin{bmatrix}
  \epsilon_i \epsilon_i \d^-_i \d^-_i x & \Gamma^+_i\d^-\dip x \\
  \eps_i\d^-_i x & x \\
  \Gamma^+_i \d^-\dip y & \Gamma^-_i \d^+\dip x \\
  y & \eps_i\d^+_iy \\
  \Gamma^-_i\d^+\dip y & \epsilon_i \epsilon_i \d^+_i \d^+_i y
\end{bmatrix} \quad \directs{i}{i+1}
$$ (Note that $\epsilon_i \epsilon_i \d^-_i \d^-_i x$ and
$\epsilon_i \epsilon_i \d^+_i \d^+_i y$ are identities for
$\circ_{i+1}$ because $\epsilon_i \epsilon_i = \epsilon_{i+1}
\epsilon_i)$.
The other case follows from the identities in Section 2.
\end{proof}

Because of Proposition~\ref{6.1}, it is convenient to regard
$\epsilon_i\d^\alpha_i$ and $\epsilon_i\d^\alpha_{i+1}$ as
generalisations of~$\psi_i$. We extend this idea to $\Psi_r$
and~$\Phi_m$, and arrive at the following definition.

\begin{definition} \label{6.2} A {\it generalised~$\psi_i$} is an
operator of the form $\psi_i$ or $\epsilon_i\d^\alpha_i$ or
$\epsilon_i\d^\alpha_{i+1}$. A {\it generalised~$\Psi_r$} is an
operator of the form $\psi'_{r-1}\psi'_{r-2}\ldots\psi'_1$, where
$\psi'_i$ is a generalised~$\psi_i$. A {\it generalised~$\Phi_m$}
is an operator of the form $\Psi'_1\Psi'_2\ldots\Psi'_m$, where
$\Psi'_r$ is a generalised~$\Psi_r$.
\end{definition}

From (2.3) and (2.5), there are results for $\epsilon_i$
and~$\d^\alpha_i$ analogous to Proposition~\ref{6.1}:
$\epsilon_i(x\circ _j y)$ is a composite of $\epsilon_i x$ and
$\epsilon_i y$; if $j=i$ then $\d^\alpha_i(x\circ _j y)$ is
$\d^\alpha_i x$ or $\d^\alpha_i y$; if $j\neq i$ then
$\d^\alpha_i(x\circ _j y)$ is a composite of $\d^\alpha_i x$ and
$\d^\alpha_i y$. From these observations and from
Proposition~\ref{6.1} we immediately get the following result.

\begin{proposition}\label{6.3} Let $\psi'_i$ be a
generalised~$\psi_i$. Then $\psi'_i(x^-\circ _j x^+)$ is naturally
equal to a composite of factors $\psi''_i x^\alpha$ with
$\psi''_i$ a generalised~$\psi_i$.

Let $\Psi'_r$ be a generalised~$\Psi_r$. Then $\Psi'_r(x^-\circ _j
x^+)$ is naturally equal to a composite of factors $\Psi''_r
x^\alpha$ with $\Psi''_r$ a generalised~$\Psi_r$.

Let $\Phi'_m$ be a generalised~$\Phi_m$. Then $\Phi'_m(x^-\circ _j
x^+)$ is naturally equal to a composite of factors $\Phi''_m
x^\alpha$ with $\Phi''_m$ a generalised~$\Phi_m$.
\end{proposition}

We will eventually express a generalised~$\Phi_n$ in terms of the
genuine folding operators~$\Phi_m$. In order to do this, we now
investigate the faces of generalised foldings.

\begin{proposition}\label{6.4} Let $\psi'_j$ be a
generalised~$\psi_j$. If $i<j$, then
$\d^\alpha_i\psi'_j=\psi''_{j-1}\d^\alpha_i$ with $\psi''_{j-1}$ a
generalised~$\psi_{j-1}$. If $i=j$, then $\d^\alpha_i\psi'_j x$ is
naturally equal to $\d^\beta_j x$ or $\d^\beta_{j+1}x$ for
some~$\beta$, or to a composite of two such factors. If $i>j$ then
$\d^\alpha_i\psi'_j=\psi''_j\d^\beta_i$ for some~$\beta$, with
$\psi''_j$ a generalised~$\psi_j$.
\end{proposition}

\begin{proof}We use relations from Section 2 and
Proposition~\ref{3.3}.

For $i<j$ we have $\d^\alpha_i\psi_j=\psi_{j-1}\d^\alpha_i$ or
$\d^\alpha_i\epsilon_j\d^\gamma_j=\epsilon_{j-1}\d^\alpha_i\d^\gamma_j
=\epsilon_{j-1}\d^\gamma_{j-1}\d^\alpha_i$ or
$\d^\alpha_i\epsilon_j\d^\gamma_{j+1}
=\epsilon_{j-1}\d^\alpha_i\d^\gamma_{j+1}
=\epsilon_{j-1}\d^\gamma_j\d^\alpha_i$.

For $i=j$ we have $\d^-_j\psi_j x=\d^-_j x\circ _j\d^+_{j+1}x$ or
$\d^+_j\psi_j x=\d^-_{j+1}x\circ _j\d^+_j x$ or
$\d^\alpha_j\epsilon_j\d^\gamma_j x=\d^\gamma_j x$ or
$\d^\alpha_j\epsilon_j\d^\gamma_{j+1}x=\d^\gamma_{j+1}x$.

For $i=j+1$ we have
$\d^\alpha_{j+1}\psi_j=\epsilon_j\d^\alpha_j\d^\alpha_{j+1}$ or
$\d^\alpha_{j+1}\epsilon_j\d^\gamma_j
=\epsilon_j\d^\alpha_j\d^\gamma_j=\epsilon_j\d^\gamma_j\d^\alpha_{j+1}$
or $\d^\alpha_{j+1}\epsilon_j\d^\gamma_{j+1}
=\epsilon_j\d^\alpha_j\d^\gamma_{j+1}$.

For $i>j+1$ we have $\d^\alpha_i\psi_j=\psi_j\d^\alpha_i$ or
$\d^\alpha_i\epsilon_j\d^\gamma_j=\epsilon_j\d^\alpha_{i-1}\d^\gamma_j
=\epsilon_j\d^\gamma_j\d^\alpha_i$ or
$\d^\alpha_i\epsilon_j\d^\gamma_{j+1}
=\epsilon_j\d^\alpha_{i-1}\d^\gamma_{j+1}
=\epsilon_j\d^\gamma_{j+1}\d^\alpha_i$.
\end{proof}

For a generalised~$\Psi_r$ we get the following results.

\begin{proposition}\label{6.5} Let $\Psi'_r$ be a
generalised~$\Psi_r$. If $i\geq r$, then
$\d^\alpha_i\Psi'_r=\Psi''_r\d^\beta_i$ for some~$\beta$, with
$\Psi''_r$ a generalised~$\Psi_r$. If $i<r$, then
$\d^\alpha_i\Psi'_r x$ is naturally equal to a composite of
factors $\Psi''_{r-1}\d^\beta_h x$ with $h\leq r$ and with
$\Psi''_{r-1}$ a generalised~$\Psi_{r-1}$.
\end{proposition}

\begin{proof}If $i\geq r$ then
$\d^\alpha_i\Psi'_r=\d^\alpha_i(\psi'_{r-1}\ldots\psi'_1)$ with
$\psi'_j$ a generalised~$\psi_j$, and the result is immediate from
Proposition~\ref{6.4}.

Now suppose that $i<r$. Then
$$\d^\alpha_i\Psi'_r x
 =\d^\alpha_i(\psi'_{r-1}\ldots\psi'_{i+1})\psi'_i\Psi'_i x,$$
with $\psi'_j$ a generalised~$\psi_j$ and with $\Psi'_i$ a
generalised~$\Psi_i$. By Proposition~\ref{6.4}
$$\d^\alpha_i\Psi'_r x
 =(\psi''_{r-2}\ldots\psi''_i)\d^\alpha_i\psi'_i\Psi'_i x$$
with $\psi''_j$ a generalised~$\psi_j$. By Propositions \ref{6.4}
and~\ref{6.3}, this is a composite of factors of the form
$$(\psi'''_{r-2}\ldots\psi'''_i)\d^\gamma_h\Psi'_i x,$$
with $i\leq h\leq i+1\leq r$ and with $\psi'''_j$ a
generalised~$\psi_j$. Since $h\geq i$, it follows from the case
already covered that the factors can be written as
$$(\psi'''_{r-2}\ldots\psi'''_i)\Psi''_i\d^\beta_h x,$$
with $\Psi''_i$ a generalised~$\Psi_i$. The factors now have the
form $\Psi''_{r-1}\d^\beta_h x$ with $\Psi''_r$ a
generalised~$\Psi_r$, as required.
\end{proof}

By iterating Proposition~\ref{6.5}, we get the following result.

\begin{proposition} \label{6.6}If $i\leq m\leq n$ and $\Psi'_r$ is a
generalised~$\Psi_r$ for $m<r\leq n$, then
$\d^\alpha_i(\Psi'_{m+1}\ldots\Psi'_n)x$ is naturally equal to a
composite of factors $(\Psi''_m\ldots\Psi''_{n-1})\d^\beta_h x$
with $h\leq n$ and with $\Psi''_r$ a generalised~$\Psi_r$.
\end{proposition}

\begin{proof}This follows from Propositions 6.5 and~\ref{6.3}.
\end{proof}

Now let $\Phi'_n$ be a generalised~$\Phi_n$; we aim to
express~$\Phi'_n$ in terms of the genuine folding
operators~$\Phi_m$. If $n=0$ or $n=1$, then necessarily
$\Phi'_n=\Phi_n$ already. In general, we use an inductive process;
the inductive step is as follows.

\begin{proposition} \label{6.7}Let $\Phi'_n$ be a generalised $\Phi_n$
which is distinct from~$\Phi_n$. Then $\Phi'_n x$ is naturally a
composite of factors $\epsilon_1\Phi'_{n-1}\d^\beta_h x$ with
$\Phi'_{n-1}$ a generalised~$\Phi_{n-1}$.
\end{proposition}

\begin{proof}By considering the first place where
$\Phi'_n$~and~$\Phi_n$ differ, we see that
$$\Phi'_n x
 =[\Phi_{m-1}(\psi_{m-1}\ldots\psi_{j+1})
 \epsilon_j][\d^\alpha_i\Psi'_j(\Psi'_{m+1}\ldots\Psi'_n)x]$$
for some $m$~and~$j$ such that $1\leq j<m\leq n$, with $i=j$ or
$i=j+1$ and with $\Psi'_r$ a generalised~$\Psi_r$. Since $j\leq
m-1$, it follows from Proposition~\ref{3.3} that
$$\Phi_{m-1}(\psi_{m-1}\ldots\psi_{j+1})\epsilon_j
 =\Phi_{m-1}\epsilon_j(\psi_{m-2}\ldots\psi_j)
 =\epsilon_1\Phi_{m-2}(\psi_{m-2}\ldots\psi_j);$$
since $j\leq i\leq m$, it follows from Propositions 6.5, 6.6
and~\ref{6.3} that
$$\d^\alpha_i\Psi'_j(\Psi'_{m+1}\ldots\Psi'_n)x$$ is a composite
of factors $\Psi''_j(\Psi''_m\ldots\Psi''_{n-1})\d^\beta_h x$
with $\Psi''_r$ a generalised~$\Psi_r$. By Proposition~\ref{6.3},
$\Phi'_n x$ is then a composite of factors of the form
$$\epsilon_1\Phi'_{m-2}(\psi'_{m-2}\ldots\psi'_j)\Psi''_j
 (\Psi''_m\ldots\Psi''_{n-1})\d^\beta_h x$$
with $\Phi'_{m-2}$ a generalised~$\Phi_{m-2}$ and with $\psi'_k$ a
generalised~$\psi_k$. These factors have the form
$\epsilon_1\Phi'_{n-1}\d^\beta_h x$ with $\Phi'_{n-1}$ a
generalised~$\Phi_{n-1}$, as required.
\end{proof}

We can now describe the interaction of~$\Phi_n$ with compositions
and face maps in general terms as follows.

\begin{proposition}\label{6.8} If a composite $x^-\circ _i x^+$ exists,
then $\Phi_n(x^-\circ _i x^+)$ is naturally equal to a composite of
factors $\epsilon_1^{n-m}\Phi_m Dx^\alpha$ with $D$ an
$(n-m)$-fold product of face operators. If $i\leq n$ then
$\d^\alpha_i\Phi_n x$ is naturally equal to a composite of factors
$\epsilon_1^{n-m-1}\Phi_m Dx$ with $D$ an $(n-m)$-fold product of
face operators.
\end{proposition}

\begin{proof}The result for $\Phi_n(x^-\circ _i x^+)$ comes from
Proposition~\ref{6.3} by iterated application of Proposition~\ref{6.7}; recall
from (2.5) that $\epsilon_1(y^-\circ _j y^+)$ is a
composite of $\epsilon_1 y^-$ and $\epsilon_1 y^+$.

Now suppose that $i\leq n$. By Proposition 3.3(iii)
$$\d^\alpha_i\Phi_n x
 =\d^\alpha_i\Phi_i (\Phi_{i+1}\ldots\Phi_n)x
 =\epsilon_1^{i-1}(\d^\alpha_1)^i(\Phi_{i+1}\ldots\Phi_n)x.$$
From Proposition~\ref{6.6}, this is a composite of factors
$\epsilon_1^{i-1}\Phi'_{n-i}D'x$ with $\Phi'_{n-i}$ a
generalised~$\Phi_{n-i}$ and with $D'$ an $i$-fold product of face
operators. By repeated application of Proposition~\ref{6.7}, there
is a further decomposition into factors $\epsilon_1^{n-m-1}\Phi_m
Dx$ with $D$ an $(n-m)$-fold product of face operators.

This completes the proof.
\end{proof}

We will now specify the composites in Proposition~\ref{6.8} more
precisely. Let $G$ be a cubical $\omega$-category, and consider
$\Phi_n(x^-\circ _i x^+)$, where $x^-\circ x^+$ is a composite
in~$G_n$. The factors $\epsilon_1^{n-m}\Phi_m Dx^\alpha$ lie in
$\Phi_n(G_n)$ (see Proposition~\ref{3.7}), and their composite can
be regarded as a composite in the $\omega$-category $\Phi_n(G_n)$
(see Theorem~\ref{3.8}). To identify the composite, we take the
universal case $$G=\lambda M(I^{i-1}\times[0,2]\times I^{n-i});$$
we may then identify $\Phi_n(G_n)$ with
$M(I^{i-1}\times[0,2]\times I^{n-i})$ by Theorem~\ref{4.1}. The
universal elements
$$x^\alpha\in[\lambda M(I^{i-1}\times[0,2]\times I^{n-i})]_n
 =\Hom[M(I^n),M(I^{i-1}\times[0,2]\times I^{n-i})]$$
are the inclusions~$\check\iota^\alpha_i$ representing
$M(I^{i-1}\times[0,2]\times I^{n-i})$ as a push-out. Evaluating
$\Phi_n(x^-\circ _i x^+)$ and the corresponding composite on~$I^n$,
and using Proposition~\ref{3.2}, we see that $\Phi_n(x^-\circ  x^+)$ gives
us $I^{i-1}\times[0,2]\times I^{n-i}$ and the factors give us
cells in $I^{i-1}\times[0,2]\times I^{n-i}$. The composite for
$\Phi_n(x^-\circ _i x^+)$ in Proposition~\ref{6.8} is an $\omega$-category
formula expressing $I^{i-1}\times[0,2]\times I^{n-i}$ as a
composite of cells. All such formulae are equivalent in all
$\omega$-categories because of the presentation of
$M(I^{i-1}\times[0,2]\times I^{n-i})$ in Theorem \ref{1.3}. The formula
uses~$\c_p$ only for $0\leq p<n$ (see Theorem~\ref{3.8}). Similarly, the
formula for $\d^\alpha_i\Phi_n x$ is an $\omega$-category formula
expressing $d^\alpha_{n-i}I^n$ as a composite of cells (see
Propositions 3.6 and~\ref{3.2}).

In order to state these results more clearly, we introduce the
following notation.

\begin{definition} \label{6.9} {\em  Let $\sigma$ be a cell in~$I^n$, and
let the dimension of~$\sigma$ be~$m$. Then   $\d_\sigma:G_m \to
G_n$ is the cubical $\omega$-category operation of the form
$\d^{\alpha(1)}_{i(1)}\ldots \d^{\alpha(n-m)}_{i(n-m)}$ such that
the underlying homomorphism
$$\check\d_\sigma\colon M(I^m)\to M(I^n)$$
sends $I^m$ to $\sigma$. }
\end{definition}

Note that $\d_\sigma$ is uniquely determined by~$\sigma$, because
of relation 2.1(i).

In this notation, we can state the following theorem.

\begin{theorem} \label{6.10} {\rm (i)} Let $f$ be a formula expressing
$I^{i-1}\times[0,2]\times I^{n-i}$ as a
$(\c_0,\ldots,\c_{n-1})$-composite of cells
$\check\iota^-_i(\sigma)$ and $\check\iota^+_i(\sigma)$, where the
$$\check\iota^\alpha_i\colon
 M(I^n)\to M(I^{i-1}\times[0,2]\times I^{n-i})$$
are the inclusions expressing $M(I^{i-1}\times[0,2]\times
I^{n-i})$ as a push-out. Let $x^-\circ _i x^+$ be a composite in a
cubical $\omega$-category. Then
$\Phi_n(x\circ _i y)$ can be got from~$f$ by replacing
$\check\iota^\alpha_i(\sigma)$ with
$\epsilon_1^{n-m}\Phi_m\d_\sigma x^\alpha$, where $m=\dim\sigma$,
and by replacing~$\c_p$ with $\circ _{n-p}$.

{\rm (ii)} Let $g$ be a formula expressing $d^\alpha_{n-i}I^n$ as
a $(\c_0,\ldots,\c_{n-2})$-composite of cells, where $1\leq i\leq
n$. In a cubical $\omega$-category,
$\d^\alpha_i\Phi_n x$ can be got from~$g$ by replacing~$\sigma$
with $\epsilon_1^{n-m-1}\Phi_m\d_\sigma x$, where $m=\dim\sigma$,
and by replacing~$\c_p$ with $\circ _{n-p}$.
\end{theorem}

\section{The natural homomorphism $B\colon G\to\lambda\gamma
G$}

Let $G$ be a cubical $\omega$-category. We
will now use Theorem~\ref{6.10} to construct a natural homomorphism
$B\colon G\to\lambda\gamma G$. Let $x$ be a member of~$G_n$. We
must define
$$B(x)\in(\lambda\gamma G)_n=\Hom[M(I^n),\gamma G].$$
Now, $M(I^n)$ is generated by the cells in~$I^n$ (see
Theorem \ref{1.3}), and $\gamma G$ is the colimit of the sequence
$$\xymatrix@1{
 \Phi_0(G_0)\ar[r]^{\epsilon_1} &
 \Phi_1(G_1)\ar[r]^{\epsilon_1} &
 \Phi_2(G_2)\ar[r] &
 \cdots}$$
(see Definition~\ref{3.9}). We can therefore define $B(x)$ by giving a
suitable value to $[B(x)](\sigma)$ for $\sigma$ a cell in~$I^n$;
these values must lie in the $\Phi_m(G_m)$, and a value
$\epsilon_1^s y$ can be identified with~$y$. The precise result is
as follows.

\begin{theorem} \label{7.1} There is a natural homomorphism $B\colon
G\to\lambda\gamma G$ for $G$ a cubical $\omega$-category given by
$$[B(x)](\sigma)=\Phi_m\d_\sigma x$$
for $\sigma$ a cell in~$I^n$, where $m=\dim\sigma$.
\end{theorem}

\begin{proof}We first show that the values prescribed for the
$[B(x)](\sigma)$ really define a homomorphism on $M(I^n)$; in
other words, we must show that they respect the relations given in
Theorem \ref{1.3}. Let $\sigma$ be an $m$-dimensional cell in~$I^n$. We
must show that $d^\alpha_m(\Phi_m\d_\sigma x)=\Phi_m\d_\sigma x$;
if $m>0$ we must also show that $d^\alpha_{m-1}(\Phi_m\d_\sigma
x)$ is the appropriate composite of the $\Phi_l\d_\tau x$, where
$\tau\subset\sigma$.

The first of these equations, $d^\alpha_m(\Phi_m\d_\sigma
x)=\Phi_m\d_\sigma x$, is an immediate consequence of Theorem~\ref{3.8}.

For the second equation, let $\sigma$ be a cell of positive
dimension~$m$. By Theorem~\ref{3.8},
$$d^\alpha_{m-1}(\Phi_m\d_\sigma x)
 =\epsilon_1\d^\alpha_1\Phi_m\d_\sigma x,$$
which may be identified with $\d^\alpha_1\Phi_m\d_\sigma x$. By
Theorem~\ref{6.10}, this is the appropriate composite of the
$\Phi_l\d_\tau x$, as required.

We have now constructed functions $B\colon G_n\to(\lambda\gamma
G)_n$, and we must show that these functions form a homomorphism
of cubical $\omega$-categories. We must
therefore show that $B(\d^\alpha_i x)=\d^\alpha_i B(x)$, that
$B(\epsilon_i x)=\epsilon_i B(x)$, that $B(x^-\circ _i
x^+)=B(x^-)\circ _i B(x^+)$, and that $B(\Gamma^{\alpha}_i
x)=\Gamma^{\alpha}_i B(x)$.

First we consider $B(\d^\alpha_i x)$, where $x\in G_n$. Let
$\sigma$ be a cell in~$I^{n-1}$ of dimension~$m$, and let
$\tau=\check\d^\alpha_i(\sigma)$. We then have
$\tau=\check\d^\alpha_i\check\d_\sigma(I^m)$, so
$\check\d_\tau=\check\d^\alpha_i\check\d_\sigma$ and
$\d_\tau=\d_\sigma\d^\alpha_i$. It follows that
$$[B(\d^\alpha_i x)](\sigma)
 =\Phi_m\d_\sigma\d^\alpha_i x
 =\Phi_m\d_\tau x$$
and
$$[\d^\alpha_i B(x)](\sigma)
 =[B(x)][\check\d^\alpha_i(\sigma)]
 =[B(x)](\tau)
 =\Phi_m\d_\tau x;$$
Therefore $[B(\d^\alpha_i x)](\sigma)=[\d^\alpha_i B(x)](\sigma)$
as required.

Next we consider $B(\epsilon_i x)$, where $x\in G_n$. Let $\sigma$
be a cell in~$I^{n+1}$ of dimension~$m$. From
Definition~\ref{2.1}, we see that $\d_\sigma\epsilon_i$ has the
form $\id\d_\tau$ or $\epsilon_j\d_\tau$. Let $l=\dim\tau$, so
that $l=m$ in the first case and $l=m-1$ in the second case. Let
$\theta\colon G_l\to G_m$ be $\id$ or $\epsilon_j\colon G_l\to
G_m$ as the case may be, and let $\check\theta\colon M(I^m)\to
M(I^l)$ be the underlying $\omega$-category homomorphism. We now
see that $\d_\sigma\epsilon_i=\theta\d_\tau$ and
$\check\epsilon_i\check\d_\sigma=\check\d_\tau\check\theta$ with
$\check\theta(I^m)=I^l$. It follows that
$$\check\epsilon_i(\sigma)
 =\check\epsilon_i\check\d_\sigma(I^m)
 =\check\d_\tau\check\theta(I^m)
 =\check\d_\tau(I^l)
 =\tau.$$
Using Theorem~\ref{5.6}, we also see that
$\Phi_m\theta=\epsilon_1^{m-l}\Phi_l$. We now get
$$[B(\epsilon_i x)](\sigma)
 =\Phi_m\d_\sigma\epsilon_i x
 =\Phi_m\theta\d_\tau x
 =\epsilon_1^{m-l}\Phi_l\d_\tau x
 =\Phi_l\d_\tau x$$
(recall that $\epsilon_1^s y$ is to be identified with~$y$) and
$$[\epsilon_i B(x)](\sigma)
 =[B(x)][\check\epsilon_i(\sigma)]
 =[B(x)](\tau)= \Phi_i\d_\tau x,$$
so that $[B(\epsilon_i x)](\sigma)=[\epsilon_i B(x)](\sigma)$ as
required.

Next we consider $B(x^-\circ _i x^+)$, where $x^-\circ _i x^+$ is
a composite in~$G_n$. Let $\sigma$ be a cell in~$I^n$ of
dimension~$m$. From Definition~\ref{2.1},
$$[B(x^-\circ _i x^+)](\sigma)=\Phi_m\d_\sigma(x^-\circ _i x^+)$$
is equal to $\Phi_m\d_\sigma x^-$ or $\Phi_m\d_\sigma x^+$ or to
$\Phi_m(\d_\sigma x^-\circ _j\d_\sigma x^+)$ for some~$j$. In any
case, using Theorem~\ref{6.10} if necessary, we see that $[B(x^-\circ _i
x^+)](\sigma)$ is a composite of factors $\Phi_l\d_\tau x^\alpha$
such that $\check\mu_i(\sigma)$ is the corresponding composite of
the $\check\iota^\alpha_i(\tau)$, where
$$\check\iota^-_i,\check\iota^+_i\colon M(I^n)\to
 M(I^{i-1}\times[0,2]\times I^{n-i})$$
are the functions expressing $M(I^{i-1}\times[0,2]\times I^{n-i})$
as a push-out. Let
$$\bigl(B(x^-),B(x^+)\bigr)\colon
 M(I^{i-1}\times[0,2]\times I^{n-i})\to\gamma G$$
be the function such that
$$\bigl(B(x^-),B(x^+)\bigr)\check\iota^\alpha_i=B(x^\alpha);$$
we see that
$$[B(x^-\circ _i x^+)](\sigma)
 =\bigl(B(x^-),B(x^+)\bigr)\check\mu_i(\sigma)
 =[B(x^-)\circ _i B(x^+)](\sigma)$$
as required.

Finally we consider $B(\Gamma^{\alpha}_i x)$, where $x\in G_n$.
Let $\sigma$ be a cell in~$I^{n+1}$ of dimension~$m$. From
Definition~\ref{2.1}, $\d_\sigma\Gamma^{\alpha}_i$ has the form
$\d_\tau$ or $\epsilon_i\d_\tau$ or $\Gamma^{\alpha}_i\d_\tau$. We
can now use the same argument as for $B(\epsilon_i x)$, noting
that $\check\Gamma^{\alpha}_i(I^m)=I^{m-1}$ and that
$\Phi_m\Gamma^{\alpha}_i=\epsilon_1\Phi_{m-1}$ by
Theorem~\ref{5.6}.

This completes the proof.
\end{proof}

\section{The natural isomorphism $B\colon G\to\lambda\gamma
G$}

In Theorem~\ref{4.1}, we have constructed a natural isomorphism
$A\colon\gamma\lambda X\to X$ for $X$ an $\omega$-category. In
Theorem~\ref{7.1} we have constructed a natural homomorphism
$B\colon G\to\lambda\gamma G$ for $G$ a cubical $\omega$-category.
We will now show that $\omega$-categories and cubical
$\omega$-categories are equivalent by showing that $B$ is an
isomorphism.

We begin with the following observation.

\begin{proposition} \label{8.1} Let $G$ be a cubical $\omega$-category.
Then $\gamma B\colon\gamma G\to\gamma\lambda\gamma G$ is an
isomorphism.
\end{proposition}

\begin{proof}Consider the composite
 $$A\circ(\gamma B)\colon\gamma G\to\gamma G.$$
By Theorem~\ref{4.1}, $A$~is an isomorphism; it therefore suffices
to show that the composite $A\circ(\gamma B)$ is the identity.
This amounts to showing that $AB(x)=x$ for $x\in\Phi_n(G_n)$. Now,
from the definitions of $A$ and~$B$, we find that
 $$AB(x)=[B(x)](I^n)=\Phi_n x;$$
since $x\in\Phi_n(G_n)$ and $\Phi_n$ is idempotent
(Proposition~\ref{3.5}), it follows that $AB(x)=x$ as required.
This completes the proof.
\end{proof}

Because of Proposition~\ref{8.1}, to show that $B$ is an
isomorphism it suffices to show that a cubical
$\omega$-category~$G$ is determined by the $\omega$-category
$\gamma G$. Because of Remark~\ref{3.10}, this is the same as
showing that $G$~is determined by the $\Phi_n(G_n)$. We will work
inductively, showing that an element~$x$ of~$G_n$ is determined by
$\Phi_n x$ and by its faces. To handle the family of faces of~$x$,
we will use the following terminology.

\begin{definition} \label{8.2}{\em  Let $G$ be a cubical $\omega$-category
and let $n$ be a positive integer. An {\it $n$-shell\/} in~$G$ an
ordered $(2n)$-tuple $$z=(z^-_1,z^+_1,\ldots,z^-_n,z^+_n)$$ of
members of~$G_{n-1}$ such that $\d^\alpha_i
z^\beta_j=\d^\beta_{j-1}z^\alpha_i$ whenever $i<j$. The set of
$n$-shells is denoted~$\Box G_{n-1}$. }
\end{definition}

\begin{remark}{\em  This construction is used in \cite[Section 5]{BH2} to
construct a coskeleton functor from $(n-1)$-truncated cubical
$\omega$-groupoids  to $n$-truncated $\omega$-groupoids determined
by $$(G_0, G_1, \ldots, G_{n-1})\mapsto (G_0, G_1, \ldots,
G_{n-1},\Box G_{n-1}),$$ and the same construction clearly works
for the category case. It follows that the folding operations are
also defined on $\Box G_{n-1}$. In the following we take a
slightly more direct route. }
\end{remark}

First, by Definition~\ref{2.1}, it is easy to check the following
result.

\begin{proposition}\label{8.3} Let $G$ be a cubical $\omega$-category
and let $n$ be a positive integer. There is a boundary map
$\d\colon G_n\to \Box G_{n-1}$ given by $$\d x=(\d^-_1 x,\d^+_1
x,\ldots,\d^-_n x,\d^+_n x).$$
\end{proposition}

Now we define folding operations on shells directly.

\begin{proposition}\label{8.4}  Let $G$ be a cubical $\omega$-category.
For $1\leq j\leq n-1$ the cubical structure of
$(G_0,\ldots,G_{n-1})$ yields a natural function $\psi_j\colon
\Box G_{n-1}\to \Box G_{n-1}$ such that $$\psi_j\d=\d\psi_j\colon
G_n\to \Box G_{n-1}.$$
\end{proposition}

\begin{proof}Let $z=(z^\alpha_i)$ be an $n$-shell. Guided by
Proposition 3.3(i), we let $\psi_i z$ be the $(2n)$-tuple
$w=(w^\alpha_i)$ such that
$$  w^\alpha_i=\begin{cases}
    \psi_{j-1}z^\alpha_i\quad & \text{for } i<j,\\
  z^-_j\circ _j z^+_{j+1}& \text{for } (\A,i)=(-,j) ,\\
   z^-_{j+1}\circ _j z^+_j  & \text{for } (\A,i)=(+,j)  ,\\
    \epsilon_j\d^\alpha_j z^\alpha_{j+1}& \text{for } i=j+1  ,\\
     \psi_j z^\alpha_i & \text{for } i>j+1.
\end{cases}
$$
From Proposition 3.1(i) and the identities in Section 2, it is
straightforward to check that $\psi_j$ is a well-defined function
from~$\Box G_{n-1}$ to itself, and it is easy to see that
$\psi_j\d=\d\psi_j$.
\end{proof}

We will now show that the $n$-dimensional elements $(n>0)$ in a
cubical $\omega$-category are determined by the lower-dimensional
elements and by the image of~$\Phi_n$.

\begin{theorem} \label{8.5}Let $G$ be a cubical $\omega$-category,
let $n$ be a positive integer, and let $\Phi_n\colon \Box G_{n-1}\to
\Box G_{n-1}$ be the function given by $$\Phi_n
 =\psi_1(\psi_2\psi_1)(\psi_3\psi_2\psi_1)\ldots
 (\psi_{n-1}\ldots\psi_1).$$
Then there is a bijection $x\mapsto(\d x,\Phi_n x)$ from~$G_n$ to
the pull-back
$$\Box G_{n-1}\times_{G_n}\Phi_n(G_n)
 =\{\,(z,y)\in \Box G_{n-1}\times\Phi_n(G_n):\Phi_n z=\d y\,\}.$$
\end{theorem}

\begin{proof}This amounts to showing that
$$\Phi_n\colon\d^{-1}(z)\to\d^{-1}(\Phi_n z)$$
is a bijection for each~$z$ in~$\Box G_{n-1}$. Since $\Phi_n$ is a
composite of operators~$\psi_j$, it suffices to show that
$$\psi_j\colon\d^{-1}(z)\to\d^{-1}(\psi_j z)$$
is a bijection for each~$z$ in~$\Box G_{n-1}$.

Given $y\in\d^{-1}(\psi_j z)$, it is straightforward to check
that there is a composite
\begin{equation*}\label{theta} \theta y
 =(\epsilon_j z^-_j\circ _{j+1}\Gamma^+_j z^+_{j+1})
 \circ _j y
 \circ _j(\Gamma^-_j z^-_{j+1}\circ _{j+1}\epsilon_j
 z^+_j)\end{equation*}
and that $\theta y\in\d^{-1}(z)$. We will carry out the proof by
showing that $\theta\psi_j x=x$ for $x\in\d^{-1}(z)$ and that
$\psi_j\theta y=y$ for $y\in\d^{-1}(\psi_j z)$.

Let $x$ be a member of $\d^{-1}(z)$. Then
$$  \theta\psi_j x
 = \begin{bmatrix}
    |& \Gamma^+_jz^-_{j+1}&  \Gamma^-_jz^-_{j+1}\\
    \eps_j z^-_j & x & \eps_jz^+_j \\
    \Gamma^+_j z^+_{j+1}&  \Gamma^-_j z^+_{j+1}& |
 \end{bmatrix} \quad \directs{j}{j+1} .
$$
The first and third rows are in the image of $\eps_{j+1}$ by (2.5)
and (2.7), so they are identities for $\circ_{j+1}$ and can
therefore be omitted. This leaves the second row in which
$\eps_jz^-_j$ and  $\eps_jz^+_j$ are identities for $\circ_j$.
It follows that $\theta \psi_j x = x$.

Now let $y$ be a member of $\d^{-1}(\psi_j z)$.
By (2.2)(vi) and (2.1)(ii),
$\eps_j\d^-_j\Gamma^+_j=\eps_j\eps_j\d^-_j=\eps_{j+1}\eps_j\d^-_j$,
so
\begin{align*}
\Gamma^+_j\d^-_{j+1}\theta y &= \Gamma^+_j z ^-_{j+1} \\ &=
\eps_j\d^-_j\Gamma^+_jz^-_{j+1}\circ_j \eps_j\d^-_j\Gamma^+_j
z^-_{j+1} \circ_j \Gamma^+_jz^-_{j+1} \\ &= \eps_{j+1}\eps_j\d^-_j
z^-_{j+1}\circ_j \eps_{j+1}\eps_j \d^-_jz^-_{j+1} \circ_j \Gamma^+_jz^-_{j+1} .\\
\intertext{Similarly} \Gamma^-_j\d^+_{j+1}\theta y
&=
\Gamma^-_jz^+_{j+1}\circ_j \eps_{j+1}\eps_j \d^+_jz^+_{j+1}\circ_j
 \eps_{j+1}\eps_j \d^+_jz^+_{j+1}. \\
 \intertext{It follows that}
 \psi_j\theta y &= \Gamma^+_j \d^-_{j+1}\theta y \circ
 _{j+1}\theta y \circ _{j+1}\Gamma^-_j\d^+_{j+1}\theta y \\
 &= \begin{bmatrix}
| & | & \Gamma^+_jz^-_{j+1} \\
\eps_jz^-_j
\circ_{j+1} \Gamma+-_jz^+_{j+1}   & y & \Gamma^-_j
z^-_{j+1}\circ_{j+1}\eps_jz^+_j\\
\Gamma^-_jz^+_{j+1} &| &|
 \end{bmatrix} \quad \directs{j}{j+1}.
\end{align*}
By (2.7) and (2.5) the first and third columns are in the image of
$\eps_j$, so they are identities for $\circ_j$ and can be omitted.
This leaves the second column so that $\psi_j\theta y =y$.

This completes the proof.
\end{proof}

From Theorem~\ref{8.5} we deduce the following result.

\begin{theorem} \label{8.6} Let $f\colon G\to H$ be a morphism of
cubical $\omega$-categories such that $\gamma
f\colon\gamma G\to\gamma H$ is an isomorphism. Then $f$ is an
isomorphism.
\end{theorem}

\begin{proof}By Remark~\ref{3.10}, $f$~induces isomorphisms from
$\Phi_n(G_n)$ to $\Phi_n(H_n)$. Since $\Phi_0$ is the identity
operation, $f$~induces a bijection from~$G_0$ to~$H_0$. By an
inductive argument using Theorem~\ref{8.5}, $f$~induces a bijection
from~$G_n$ to~$H_n$ for all~$n$. Therefore $f$ is an isomorphism.
\end{proof}

It follows from Proposition~\ref{8.1} and Theorem~\ref{8.6} that
$B\colon G\to\lambda\gamma G$ is a natural isomorphism for cubical
$\omega$-categories~$G$. From Theorem~\ref{4.1},
$A\colon\gamma\lambda X\to X$ is a natural isomorphism for
$\omega$-categories~$X$. We draw the following conclusion.

\begin{theorem} \label{8.7} The categories of $\omega$-categories and
of cubical $\omega$-categories are equivalent
under the functors $\lambda$ and~$\gamma$.
\end{theorem}

\section{Thin elements and commutative shells  in a  cubical $\omega$-category}

In this section we use the equivalence of Theorem 8.8 to clarify
two concepts in the theory of cubical $\omega$-categories: thin
elements and commutative shells. Thin elements (sometimes called
hollow elements) were introduced in the thesis of K. Dakin
\cite{Da}, and were developed in the cubical $\omega$-groupoid
context by Brown and Higgins in \cite{BH1,BH2,BH2a}. They are used
by Ashley in \cite{Ash} and by Street in \cite{Street1}. In the
cubical nerve of an $\omega$-category they arise as follows.

Throughout this section, let $G$ be a cubical $\omega$-category.
Whenever convenient, we will identify $G$ with the nerve of
$\gamma G$; in other words, each element $x$ of $G_n$ is
identified with a homomorphism $x: M(I^n) \to \gamma G$.

First we deal with thin elements. Intuitively, an element is thin
if its real dimension is less than its apparent dimension. In the
nerve of an $\omega$-category we can make this precise as follows.

\begin{definition} {\em Let $x$ be a member of $G_n$. Then $x$ is
{\em thin} if $$\mathrm{dim } \; x(I^n) < n.$$}
\end{definition}

Given an element $x$ of $G_n$, we can identify $x(I^n)$ with
$\Phi_n x$ by Theorem 7.1. By Remark 3.10, dim $\Phi_n x <
n$ if and only if $\Phi_n x$ is in the image of $\epsilon_1$. We
therefore have the following characterisation.

\begin{proposition} Let $x$ be a member of $G_n$. Then $x$ is thin
if and only if $\Phi_n x$ is in the image of $\epsilon_1$.
\end{proposition}

There is also a less obvious characterisation in more
elementary cubical terms: the thin elements of $G_n$ are
those generated by the $G_m$ with $m<n$. The precise
statement is as follows.

\begin{theorem} Let $x$ be a member of $G_n$. Then $x$ is thin if
and only if it is a composite of elements of the forms
$\epsilon_i y$ and $\Gamma^\alpha_j z$ for various values of
$i, j, \alpha,\; y,\;z$.
\end{theorem}
\begin{proof} Suppose that $x$ is a composite of elements of the
forms $\epsilon_i y$ and $\Gamma^\alpha_j z$. Then $\Phi_n x$
is in the image of $\epsilon_1$ by Theorems 5.6 and 6.10,
so x is thin by Proposition 9.2.

Conversely, suppose that $x$ is thin. It follows from the
proof of Theorem 8.5 that $x$ is a composite of $\Phi_n x$
with elements of the forms $\epsilon_i y$ and
$\Gamma^\alpha_j z$. By Proposition 9.2, $\Phi_n x$ is in the
image of $\epsilon_1$, so $x$ is itself a composite of elements
of the forms $\epsilon_i y$ and $\Gamma^\alpha_j z$.
\end{proof}

Next we deal with commutative shells. There is an
obvious concept of commutative square, or commutative
2-shell; we want commutative $n$-shells for arbitrary
positive $n$. Now an $n$-shell $z$ in $G$ can be identified with a
homomorphism
$z: M(d^-_{n-1} I^n \cup d^+_{n-1} I^n) \to \gamma G$,
and we must obviously define a commutative $n$-shell as
follows.

\begin{definition} {\em For $n>0$ an $n$-shell $z$ in $G$ is {\em commutative} if
$$z(d^-_{n-1} I^n) = z(d^+_{n-1} I^n).$$ }
\end{definition}

By Theorem 6.10(ii), if $z$ is an $n$-shell with $n>0$ then
$z(d^\alpha_{n-1})$ can be identified with $(\Phi_n z)^\alpha_1$,
the $(\alpha,1)$ face of the $n$-shell $\Phi_n z$ as in theorem
8.6. We can therefore describe commutative $n$-shells in cubical
terms as follows.

\begin{proposition} For $n>0$ an $n$-shell $z$ in $G$ is commutative if
and only if $$(\Phi_n z)^-_1 = (\Phi_n z)^+_1.$$
\end{proposition}
\section{Monoidal closed structures}
In \cite{AAS} Al-Agl and Steiner constructed a  monoidal closed
structure on the category $\omega\dash\cat^\glob$ of (globular)
$\omega$-categories by using a cubical description of that
category. Now that we have a more explicit cubical description we
can give a more explicit description of the monoidal  closed
structure; we modify the construction which is  given by Brown and
Higgins in \cite{BHtens} for the case of a single connection and
for groupoids rather than categories. Following the method there,
we first  define the closed structure on the category
$\omega\dash\cat^\Box$ of cubical $\omega$-categories using a
notion of $n$-fold left homotopy which we outline below, and then
obtain the tensor product as the adjoint to the closed structure.
This gives:

\begin{theorem}The category $\omega\dash\cat^\Box$ admits a monoidal closed structure
with an adjoint relationship $$ \omega\dash\cat^\Box( G \otimes H,
K) \cong \omega\dash\cat^\Box(G, \omega\dash\Cat^\Box(H,K)) $$ in
which $\omega\dash\Cat^\Box(H,K)_0$ is the set of morphisms $H\to
K$, while for $ n \ge 1$ $\omega\dash\Cat^\Box(H,K)_n$ is the set
of $n$-fold left homotopies $H \to K$.
\end{theorem}

The proof is given below.

Because of the equivalence between $\omega\dash\cat^\Box$ and the
category $\omega\dash\cat^\glob$ of $\omega$-categories we then
have:

\begin{corollary}\label{globmonoid} The category $\omega\dash\cat^\glob$ admits a monoidal
closed structure with an adjoint relationship

$$ \omega\dash\cat^\glob( X \otimes Y, Z) \cong
\omega\dash\cat^\glob(X, \omega\dash\Cat^\glob(Y,Z)) $$
 in which    $\omega^\glob\dash\Cat(Y,Z)_0$ is
the set of morphisms $Y\to Z$,  while  for $n \ge 1$
$\omega\dash\Cat^\glob(Y,Z)_n$ is the set of $n$-fold left
homotopies $Y \to Z$ corresponding to the cubical homotopies.
\end{corollary}

The tensor product in Corollary \ref{globmonoid}  is an extension
of the tensor product in Theorem \ref{1.8}.

Note that by Remark 3.11, the set $\omega\dash\Cat^\glob(Y,Z)_n$
of globular $n$-fold left homotopies may be thought of as an
explicitly described subset of the set of cubical $n$-fold left
homotopies $\lambda Y \to \lambda Z$.   Because of the
complications of the folding operations, explicit descriptions of
the globular monoidal closed structure are not so easy, but have
been partly accomplished by Steiner in \cite{Ste-past}. See also
Crans \cite{Cr}.

We now give details of these cubical constructions, following
directly the methods of \cite{BHtens}.

Let $H$ be a \omcat y and $n$ be a non-negative integer. We can
construct a \omcat y $P^n H$ called the {\em $n$-fold (left) path
\omcat y\/} of~$H$ as follows: $(P^n H)_r=H_{n+r}$; the operations
$\d^\A_i$, $\eps_i$, $\Gamma^\A_i$ and~$\circ_i$ of $P^n H$ are
the operations $\d^\A_{n+i}$, $\eps_{n+i}$, $\Gamma^\A_{n+i}$
and~$\circ_{n+i}$ of~$H$. The operations $\d^\A_1, \ldots,
\d^\A_n$ not used in $P^n H$ give us morphisms of \omcat ies from
$P^n H$ to $P^{n-1} H$, etc., and we get an internal \omcat y
$$\P H = (H, P^1 H, P^2 H, \ldots \,)$$
in the category \ocat.

For any \omcat ies $G,H$ we now define
$$ \oCatm(G,H)=\ocatm(G,\P H);$$
that is, $\oCatm_m(G,H)= \ocatm(G,P^mH)$, and the \omcat y
structure on $\oCatm_m(G,H)$ is induced by the internal \omcat y
structure on $\P H$. Ultimately, this means that the operations
$\d^\A_i$, etc.\ on $\oCatm_m(G,H)$ are induced by the similarly
numbered operations on~$H$. In dimension $0$, $\oCatm(G,H)$
consists of all morphisms $G \to H$, while in dimension $n$ it
consists of {\em $n$-fold (left) homotopies} $G \to H$. We make
\oCat $(G,H)$ a functor in $G$ and $H$ (contravariant in $G$) in
the obvious way.

The definition of tensor product of \omcat ies is harder. We
require that $-\otimes G$ be left adjoint to \oCat$(G,-)$ as a
functor from \ocat \ to \ocat, and this determines $\otimes $ up
to natural isomorphism.   Its existence, that is, the
representability of the functor $\ocatm(F,\oCatm(G,-))$, can be
asserted on general grounds. Indeed, \ocat\  is an equationally
defined category of many sorted algebras in which the domains of
the operations are defined by finite limit diagrams, and general
theorems on such algebraic categories imply that \ocat\ is
complete and cocomplete.

We can also specify the tensor product \omcat y by a {\em
presentation\/}; that is, we give a set of generators in each
dimension and a set of relations of the form $ u=v$, where $u,v$
are well formed formulae of the same dimension made from
generators and the operators $\d^\A_i, \eps_i,
\Gamma^\A_i,\circ_i$. This is analogous to the standard tensor
product of modules over a ring, and the universal property of the
presentation gives the required adjointness.

The details are as follows.
\begin{definition}  {\em
Let $F,G$ be \omcat ies. Then $F \otimes G$ is the \omcat y
generated by elements in dimension $ n \ge 0$ of the form $ x
\otimes y$ where $ x \in F_p, \; y \in G_q $ and $ p+q =n$,
subject to the following defining relations (plus, of course, the
laws for \omcat ies):
\begin{enumerate}[(i)]
\item $\d^\A_i(x \otimes y) =
\begin{cases}(\d^\A_i x) \otimes y & \text{if } 1 \le i \le p,\\
              x \otimes (\d^\A_{i-p}y)& \text{if } p+ 1 \le i \le n;
\end{cases}$
 \item $\eps_i(x \otimes y) =
\begin{cases}(\eps_i x) \otimes y & \text{if } 1 \le i \le p+1,\\
              x \otimes (\eps_{i-p}y)& \text{if } p+ 1 \le i \le n+1;
\end{cases}$
\item $\Gamma^\A_i(x \otimes y) =
\begin{cases}(\Gamma^\A_i x) \otimes y & \text{if } 1 \le i \le p,\\
              x \otimes (\Gamma^\A_{i-p}y)& \text{if } p+ 1 \le i \le n;
\end{cases}$
\item $(x\circ_i x') \otimes y =
(x \otimes y)\circ_i (x' \otimes y)$ if $ 1 \le i \le p,$  and $
x \circ_i x'$ is defined in  $ F; $
\item $x \otimes (y\circ_j y') =
(x \otimes y)\circ_{p+j} (x \otimes y') $
 if  $1 \le j \le q,$ and  $ y \circ_j y'$ is
defined in $ G; $
\end{enumerate}     }
\end{definition}
We note that the relation
\begin{enumerate}
  \item[(vi)] $ (\eps_{p+1} x) \otimes y = x \otimes (\eps _1 y)$
  \end{enumerate}
follows from (ii).

An alternative way of stating this definition is to define a {\em
bimorphism} $(F,G) \to A$, where $F,G,A$ are \omcat ies, to be a
family of maps $F_p \times G_q \to A_{p+q} \; (p,q \ge 0)$,
denoted by $(x,y) \mapsto \chi (x,y)$ such that
\begin{enumerate}[(a)]
  \item for each $ x \in F_p$, the map $ y \mapsto \chi(x,y)$ is a
  morphism of \omcat ies $G \to P^pA$;
  \item for each $g \in G_q$ the map   $ x \mapsto \chi(x,y)$ is a
  morphism of \omcat ies $F \to TP^qTA$,
\end{enumerate}
 where  the \omcat y $TX$ has the same elements as
 $X$ but its cubical operations, connections and compositions are
 numbered in reverse order. The \omcat y $F \otimes G$ is now
 defined up to natural isomorphisms by the two properties
\begin{enumerate}[(i)]
  \item the map $(x,y) \mapsto x \otimes y$ is a bimorphism $(F,G)
  \to F \otimes G$;
  \item every bimorphism $(F,G) \to A$ is uniquely of the form
  $(x,y) \mapsto \sigma (x \otimes y)$ where $ \sigma : F \otimes
  G \to A$ is a morphism of \omcat ies.
\end{enumerate}

In the definition of a bimorphism $(F,G) \to A$, condition (a)
gives maps $F_p \to \oCatm_p(G,A)$ for each $p$, and condition (b)
states that these combine to give a morphism of \omcat ies $F \to
\oCatm(G,A)$. This observation yields a natural bijection between
bimorphisms $(F,G) \to A$ and morphisms $F \to \oCatm(G,A)$. Since
we also have a natural bijection between bimorphisms $(F,G) \to A$
and morphisms $F \otimes G \to A$, we have
\begin{proposition}
The functor $- \otimes G$ is left adjoint to the functor
$\oCatm(G,-)$ from $\ocatm$ to $\ocatm$. \hfill $\Box$
\end{proposition}

\begin{proposition} \label{10.4}
  For \omcat ies $F,G,H$, there are natural isomorphisms of \omcat
  ies
\begin{enumerate}
  \item[\rm{(i)}] $(F \otimes G) \otimes H \cong F \otimes (G \otimes H) $,
  and
  \item[\rm{(ii)}] $\oCatm(F \otimes G,H) \cong \oCatm (F, \oCatm (G,H))$
\end{enumerate}
giving $\ocatm$ the structure of a monoidal closed category.
\end{proposition}
\begin{proof}
(ii) In dimension $r$ there is by adjointness a natural bijection
\begin{align*}
\oCatm_r(F \otimes G,H) & = \ocatm(F \otimes G, P^rH) \\
                        & \cong \ocatm(F, \oCatm(G,P^rH)) \\
                        &= \ocatm(F, P^r(\oCatm(G,H)))\\
                        &= \oCatm_r(F,\oCatm(G,H)) .
\end{align*}
These bijections combine to form the natural isomorphism (ii) of
\omcat ies because, on both sides, the \omcat y structures are
induced by the corresponding operators $\d^\A_i,\eps_j$, etc. in
$H$.

(i) This isomorphism may be proved directly, or, as is well known,
be deduced from the axioms for a monoidal closed category.
\end{proof}

We can also relate the construction to the category of cubical
sets, which we denote $\cub$. The underlying cubical set functor
$U: \ocatm \to \cub$ has a left adjoint $\sigma : \cub \to
\ocatm$, and we call $\sigma(K)$ the {\em free \omcat y on the
cubical set } $K$. The category $\cub$ has a monoidal closed
structure in the same way as $\ocatm$ (see \cite{BHtens}); the
internal hom $\Cub$ is given by $\Cub(L,M)_r=\cub(L,P^rM)$ where
$P^r$ is now the $n$-fold path functor on cubical sets. We have
the following results.

\begin{proposition}
For a cubical set $L$ and \omcat y $G$, there is a natural
isomorphism of cubical sets
$$U(\oCatm(\sigma(L),G)) \cong \Cub(L,UG).$$
\end{proposition} \begin{proof}
The functor $\sigma: \cub \to \ocatm$ is left adjoint to $U :
\ocatm \to \cub$, and this is what the proposition says in
dimension 0. In dimension $r$ we have a natural bijection
 \begin{align*}
\oCatm_r(\sigma(L),G) &= \ocatm(\sigma(L),P^rG) \\
                      & \cong \cub(L,UP^rG)\\
                      &=\Cub_r(L,UG)
   \end{align*}
and these bijections are compatible with the cubical operators.
\end{proof}
\begin{proposition}
If $K,L$ are cubical sets, there is a natural isomorphism of
\omcat ies
$$  \sigma  (K) \otimes \sigma (L) \cong  \sigma(K \otimes L).$$
\end{proposition}
\begin{proof}
For any \omcat y $G$, there are natural isomorphisms of cubical
sets
\begin{align*}
U(\oCatm(\sigma(K) \otimes \sigma (L),G)) &\cong
U(\oCatm(\sigma(K),\oCatm(\sigma(L),G))) \\
   &\cong \Cub(K,U(\oCatm(\sigma(L),G))) \\
   &\cong  \Cub(K,\Cub(L,UG))\\
   &\cong \Cub(K \otimes L,UG)\\
   &\cong U(\oCatm(\sigma(K \otimes L),G).\\
\intertext{The proposition follows from the information in
dimension $0$, namely} \ocatm( \sigma(K \otimes L),G) &\cong
\ocatm(\sigma(K) \otimes \sigma (L),G)) .
\end{align*}\end{proof}
The $\omega$-categories $M(I^n)$ of Section~1 can be fitted into
this framework if one regards them as \omcat ies. Indeed, as a
\omcat y, $M(I^n)$ is freely generated by one element in
dimension~$n$; therefore $M(I^n)=\sigma(\bI^n)$ where $\bI^n$ is
the cubical set freely generated by one element in dimension~$n$.
Calculations with cubical sets show that $\bI^m \otimes \bI^n
\cong \bI^{m+n}$, and we get the following result.

\begin{corollary}
The are natural isomorphisms of \omcat ies
\begin{equation}M(I^m) \otimes M(I^n) \cong M(I^{m+n}). \tag*{$\Box$}  \end{equation}

\end{corollary}
\begin{proposition}

\begin{enumerate}
  \item[\rm(i)] $M(I^n) \otimes -$ is left adjoint to $P^n: \ocatm
  \to \ocatm$.
  \item[\rm(ii)] $-\otimes M(I^n)$ is left adjoint to
  $\oCatm(M(I^n),-)$.
  \item[\rm(iii)] $\oCatm(M(I^n),-)$ is naturally isomorphic to
  $TP^nT$.
\end{enumerate}
\end{proposition}
\begin{proof}
(i) There are natural bijections
\begin{align*}
\ocatm(M(I^n) \otimes H,K) & \cong \ocatm(M(I^n),\oCatm(H,K))
\\ & \cong \oCatm_n(H,K) \\ &=\ocatm(H,P^nK).
\end{align*}

(ii) This is a special case of Proposition \ref{10.4}.

(iii) It follows from (i) that $TP^nT: \ocatm \to \ocatm$ has left
adjoint $T(M(I^n) \otimes T(-)) \cong - \otimes TM(I^n)$. But the
obvious isomorphism $T \bI \to \bI$ induces an isomorphism
$TM(I^n) \cong M(I^n)$, so $-\otimes TM(I^n)$ is naturally
isomorphic to $- \otimes M(I^n)$. The result now follows from
(ii).
\end{proof}

The free \omcat y on a cubical set is important in applications to
concurrency theory. The data for a concurrent process can be given
as a cubical set $K$, and the evolution of the data can be
reasonably described by the free \omcat y $\sigma(K)$; indeed,
$\sigma(K)$ is the higher-dimensional analogue of the path
category on a directed graph. The idea is pursued by Gaucher in
\cite{Gauch}.

\section*{Acknowledgements} We would like to thank the Department of
Mathematical Sciences at Aalborg University for supporting Brown
at the June, 1999,  GETCO Workshop; the London Mathematical
Society for support of a small workshop on multiple categories and
concurrency in September, 1999, at Bangor; and a referee for
helpful comments. This paper is dedicated to Philip Higgins in
recognition of the power of his insights into `the algebra of
cubes' and in thanks for a long and happy collaboration with the
second author.

\end{document}